\documentclass[11pt]{amsart}

\usepackage[backref,bookmarks=false]{hyperref}

\usepackage{color}

\numberwithin{equation}{section}

\newcommand{\ben}{\begin{enumerate}}
\newcommand{\een}{\end{enumerate}}

\newcommand{\bea}{\begin{eqnarray}}
\newcommand{\ba}{\begin{array}}
\newcommand{\bean}{\begin{eqnarray*}}
\newcommand{\ea}{\end{array}}
\newcommand{\eea}{\end{eqnarray}}
\newcommand{\eean}{\end{eqnarray*}}
\newcommand{\beq}{\begin{equation}}
\newcommand{\eeq}{\end{equation}}
\newcommand{\bthm}{\begin{thm}}
\newcommand{\ethm}{\end{thm}}
\newcommand{\blem}{\begin{lem}}
\newcommand{\elem}{\end{lem}}
\newcommand{\bprop}{\begin{prop}}
\newcommand{\eprop}{\end{prop}}
\newcommand{\bcor}{\begin{cor}}
\newcommand{\ecor}{\end{cor}}
\newcommand{\bdfn}{\begin{dfn}}
\newcommand{\edfn}{\end{dfn}}
\newcommand{\brem}{\begin{rem}}
\newcommand{\erem}{\end{rem}}
\newcommand{\bpf}{\begin{proof}}
\newcommand{\epf}{\end{proof}}
\newcommand{\bfact}{\begin{fact}}
\newcommand{\efact}{\end{fact}}
\newcommand{\bobs}{\begin{obs}}
\newcommand{\eobs}{\end{obs}}
\newcommand{\bexam}{\begin{exam}}
\newcommand{\eexam}{\end{exam}}
\newcommand{\bclaim}{\begin{claim}}
\newcommand{\eclaim}{\end{claim}}

\newtheorem{thm}{Theorem}[section]
\newtheorem{prop}[thm]{Proposition}
\newtheorem{lem}[thm]{Lemma}

\newtheorem{cor}[thm]{Corollary}
\newtheorem{dfn}[thm]{Definition}
\newtheorem{rem}[thm]{Remark}
\newtheorem{fact}[thm]{Fact}
\newtheorem{claim}[thm]{Claim}
\newtheorem{obs}[thm]{Observation}
\newtheorem{exam}[thm]{Example}

\newtheorem*{condition'}{Condition 2'}

 \newtheoremstyle{claimstyle}%
   {}
   {}
   {\normalfont}
   {}
   {\itshape}
   {.}
   { }
   {\thmnote{#3}}

\theoremstyle{claimstyle}

\alph{enumii} \roman{enumiii}

\def\cA{\mathcal A}             \def\cB{\mathcal B}       \def\cC{\mathcal C}
\def\cH{\mathcal H}             \def\cF{\mathfrak F}       
\def\cL{{\mathcal L}}           \def\cM{\mathcal M}        \def\cP{{\mathcal P}}

           \def\cK{\mathcal K}

\def\N{{\mathbb N}}                \def\Z{{\mathbb Z}}      \def\R{{\mathbb R}}

\def\F{{\mathcal F}}

\def\a{\alpha}                \def\b{\beta}             
               \def\e{\varepsilon}           
\def\g{\gamma}                \def\Ga{\Gamma}           \def\l{\lambda} \def\la{\lambda}
\def\La{\Lambda}              \def\om{\omega}           \def\Om{\Omega}
\def\Sg{\Sigma}               \def\sg{\sigma}
               \def\th{\theta}           \def\vth{\vartheta}
\def\ka{\kappa}

\newcommand{\ep}{\varepsilon}
\newcommand{\ph}{\varphi}
\newcommand{\al}{\alpha}
\newcommand{\ga}{\gamma}

\def\1{1\!\!1}

\def\and{\text{ and }}

\def\dist{\text{{\rm dist}}}

              \def\bu{\bigcup}
\def\({\bigl(}                \def\){\bigr)}
\def\lt{\left}                \def\rt{\right}

\def\ld{\ldots}                        \def\^{\tilde}

\def\es{\emptyset}            \def\sms{\setminus}
\def\sbt{\subset}

\def\sp{\medskip}             \def\fr{\noindent}        
            
\def\ess{{\rm ess}}
\def\om{\omega}

\def\${$ \displaystyle }

\def\den{\rho}
\def\shift{\theta}

\newcommand{\pf}{{\mathcal{L}}}
\newcommand{\pfx}{{\mathcal{L}}_x}
\newcommand{\npf}{\mathcal{\hat L}}
\newcommand{\npfx}{{ \mathcal{\hat L}}_x}

\begin{document}

\title[Random Shifts]{Countable Alphabet Random Subhifts of finite type with weakly positive transfer operator}


\author{Volker Mayer}
\address{Universit\'e de Lille I, UFR de
  Math\'ematiques, UMR 8524 du CNRS, 59655 Villeneuve d'Ascq Cedex,
  France} \email{volker.mayer@math.univ-lille1.fr \newline
  \hspace*{0.42cm} \it Web: \rm math.univ-lille1.fr/$\sim$mayer}

\author{Mariusz Urba\'nski}
\address{Department of
  Mathematics, University of North Texas, Denton, TX 76203-1430, USA}
\email{urbanski@unt.edu \newline \hspace*{0.42cm} \it Web: \rm
  www.math.unt.edu/$\sim$urbanski}

\date{\today} \subjclass{111}

\thanks{The second-named author was supported in part by the NSF Grant DMS 1361677.}

\maketitle
\begin{abstract}
We deal with countable alphabet locally compact random subshifts of finite type (the latter merely meaning that the symbol space is generated by an incidence matrix) under the absence of Big Images Property and under the absence of uniform positivity of the transfer operator. 
We first establish the existence of random conformal measures along with good bounds for the iterates of the Perron-Frobenius operator. 
Then, using the technique of positive cones and proving a version of Bowen's type contraction (see \cite{Bow75}), we also establish a fairly complete thermodynamical formalism. This means that we prove the existence and uniqueness of fiberwise invariant measures (giving rise to a global invariant measure) equivalent to the fiberwise conformal measures. Furthermore, we establish the existence of a spectral gap for the transfer operators, which in the random context precisely means the
exponential rate of convergence of the normalized iterated transfer operator. This latter property in a relatively straightforward way entails the exponential decay of correlations and the Central Limit Theorem. 
\end{abstract}

\section{Introduction}
The thermodynamic formalism of finite and countable topological Markov shifts (subshifts of finite type) with H\"older continuous potentials is by now quite well understood. The case of finite alphabet was settled already in nineteen seventies, primarily due to the work of R. Bowen, \cite{Bow75}, and Ruelle, \cite{Rue78}. For the case of countable infinite alphabet, the existence and uniqueness of conformal measures and invariant Gibbs states for finitely irreducible shifts was established in \cite{MauldinUrb01}. The necessity of finite irreducibility for the existence of invariant Gibbs states was shown in \cite{Sarig03}. The spectral gap of the corresponding Perron-Frobenius (transfer) operator, and resulting from it exponential decay of correlations, the Central Limit Theorem and the Law of Iterated logarithm were established in \cite{MauldinUrb03}. As a matter of fact, this book contains a systematic exposition of the theory of deterministic subshifts of finite type. The work \cite{SarigLN09} is also a good reference with a different viewpoint based on the concept of so called recurrent potentials.

\sp In the case of random dynamics, Bogensch\"utz and Gundlach \cite{BogGun95} and also Kifer in \cite{Kif08} have considered the finite shift case; see also \cite{MSU} for a general and systematic treatment of random expanding maps. Then Denker, Kifer and Stadlbauer in \cite{DKS08},
Stadlbauer \cite{Stadl10, Stadl13} and \cite{RoyUrb2011} dealt with the case of random countable topological Markov shifts. All these papers assumed some version of finite irreducibility or, to put it in different words, Big Images Property. In all these cases the Perron-Frobenius evaluated at the function identically equal to one was strictly positive.

\sp The present paper lies on just the opposite spectrum. It is devoted to a class of random countable subshifts of finite type for which the Big Images Property does always fail. Also, strict positivity of the Perron-Frobenius does always fail, see Property (B) in Section~\ref{Section_TO}. This happens however with some control as can be seen from properties (A) and (C) of the same section. Our hypotheses however entail local compactness of the fiberwise symbol spaces. Such random shifts occur naturally in various situations and they are closely related to stochastic processes in random environments (see \cite{Kif98,ANV}). We describe some examples in Section \ref{sec examples}.

\sp In the present paper, utilizing the concept of narrow topology and Prokhorov's Compactness Theorem, we first establish the existence of conformal measures along with good bounds for the iterates of the Perron-Frobenius operator. 
Then, using the technique of positive cones and proving a version of Bowen's type contraction (see \cite{Bow75}), by developing in the present setting our approach from a non-Markovian context of \cite{MyUrb2014}, we also establish a fairly complete thermodynamical formalism. First of all, we prove the existence and uniqueness of fiberwise invariant measures (giving rise to a global invariant measure) equivalent to the fiberwise conformal measures. Furthermore, we establish the existence of a spectral gap for the transfer operators, which in the random context precisely means the
exponential rate of convergence of iterates of the normalized transfer operator. This latter property in a relatively straightforward way entails the exponential decay of correlations and the Central Limit Theorem

\sp We would like to stress that the hypotheses under which we work, apart from topological mixing of the shift map, are of ``the first level'' type, meaning that these do not involve iterates of the random shift map, Birkhoff's sums of potentials, or iterates of the transfer operator. These involve the transfer operator itself but only very mildly as condition (C) of Section~\ref{Section_TO}. 

\

{\bf Acknowledgement.} We wish to thank the anonymous referee for valuable remarks and suggestions which improved the final exposition of our paper. We particularly thank him or her for bringing up to our attention the distinction between annealed and quenched results.

\section{Preliminaries; Topological Aspects}

Let $(X, \cF, m)$ be a complete probability space with $\cF$, a $\sg$-algebra of subsets of $X$ and $m$ a complete probability measure on $\cF$. Let $\th:X\to X$ be an invertible measurable transformation preserving measure $m$. Let $E$ be a countable infinite set. For the ease of exposition and without loss of generality we assume throughout that
$$
E=\N=\{0,1,2,\ld\}
$$
is the set of all non-negative integers. 
Assume that 
$$
X\ni x\mapsto A(x):E\times E\to\{0,1\}
$$ 
is a measurable map from $X$ into $0,1$-matrices on $E$, commonly called incidence matrices. These matrices define symbol spaces as follows. For every $x\in X$ let
\beq\label{1.2}
E_x^\infty=E_x^\infty(A):=\big\{(\om_n)_{n=0}^\infty\in E^\N: A_{\om_i\om_{i+1}}(\th ^i (x))=1 \; \; \text{for every } i\geq 0\big\}.
\eeq
More generally, for every set $F\subset E$, given  $n\in\N\cup\{\infty\}$, it is useful to introduce the sets 
\beq\label{2.2}
F_x^n=F_x^n(A)=\big\{\om_0\om_1\ld\om_n\in F^{n+1} \; : \; A_{\om_i\om_{i+1}}(\th ^i (x))=1  \  \text{ for all } \ 0\leq i \le n-1\big\}.
\eeq
These are words of length $n+1$ over the alphabet $F$ and they will be called admissible by the matrix $A$, also $A$-admissible, or just admissible.
If $m\leq n$ and if $\om\in F_x^n$ then 
$$
\om|_m:=\om_0\om_1...\,\om_m\in F_x^m
$$
is called the truncation, or restriction, of $\om$ to $m$. The cylinder generated by a word $\om\in E_x^n$ is defined to be
$$
[\om]_x := \{ \tau\in E_x^\infty \; : \;\; \tau|_n=\om\}\,.
$$
A straightforward extension of the notion of cylinder is this. For any $F\subset E$, 
$$
[F]_x:=\bigcup_{j\in F} [j]_x,
$$
The symbol spaces $E_x^\infty$, $x\in X$, are naturally endowed with the subspace topology inherited from the product (Tichonov) topology on the Cartesian product $E^\N$. This latter topology, and all respective subspace topologies are generated by many natural metrics. One of them, the one we will work with, is this:
$$ 
d(\om,\tau)= 0 \;\;\text{if}\; \; \om=\tau \;\;\text{and}\; \;
d(\om,\tau)=\exp\Big(-\min \{n\;  ; \;\om_n\neq \tau_n\} \Big)\;\; \text{if}\;\; \om\neq \tau\,.
$$
The Cartesian product $X \times E^{\N}$ is further endowed with a natural measurable structure, i.e. a $\sg$-algebra. This is the $\sg$-algebra generated by the sets of the form $F\times B$, where $F$ is a measurable subset of $X$ and $B$ is a Borel subset of $E^\N$. This $\sg$-algebra is denoted by $\cF\otimes\cB$. Let
\beq\label{2.3} 
\Om:=\bigcup_{x\in X} \{x\}\times E_x^\infty \subset X \times E^{\N}
\eeq
By $\cF\otimes|_\Om\cB$ we denote the restriction of $\cF\otimes\cB$ to $\Om$, i.e. the $\sg$-algebra consisting of the sets $Y\cap\Om$, $Y\in\cF\otimes\cB$. By $\cB_x$, $x\in X$, we denote the Borel $\sg$-algebra of Borel subsets of $E_x^\infty$. The symbol $\sg$ stands for the standard shift map from $E^\N$ to $E^\N$ defined by the formula
$$
\sg\((\om_n)_{n=0}^\infty\)=(\om_{n+1})_{n=0}^\infty.
$$
The spaces $E_x^\infty$, $x\in X$, are invariant under the shift map $\sg:E^\N\to E^\N$ in the sense that
$$
\sg(E_x^\infty)\sbt E_{\th(x)}^\infty
$$
for all $x\in X$. We then set
$$
\sg_x:=:\sg{\big|}_{E_x^\infty}:E_x^\infty\to E_{\th(x)}^\infty.
$$
We are interested in the (fiberwise) dynamics, i. e. in the compositions (frequently referred to as iterates):
$$
\sg_x^n 
:= \sg_{\th ^{n-1}(x)}\circ...\circ \sg_x:E_x^\infty \to E_{\th^n(x)}^\infty ,  \  \ n\geq 0.
$$
Its global version, the random shift map 
$
\hat\sg_A:\Om\to\Om,
$
defined by the formula
\beq\label{1.21}
(x,\om) \longmapsto (\th (x) , \sg_x (\om)), 
\eeq
is a skew product map with base $X$ and fibers $E_x^\infty$. It will be frequently denoted simply by $\hat\sg$. We assume this map to be mixing. The precise definition is the following.

\bdfn\label{2.4}
The random shift map $\hat\sg:\Om\to\Om$ is called topologically mixing if for all letters $a,b\in E$ there exists $N=N_{a,b}\ge 0$ such that for every $n\geq N$
and all $x\in X$ there exists $\om\in E_x^n$ such that the word $a\om b$ is  
$A$-admissible. 
\edfn

\fr A more adequate name for this concept would be uniform topological mixing, but we stick to the shorter one since there will be considered in this paper no more general, or weaker, concepts of mixing.

\brem\label{2.001} Notice that mixing in particular implies that no cylinders $[i]_x$ and no sets $\sg_x^{-1}\([i]_{\th(x)}\)$, $i\in E$, $x\in X$, are empty. A reformulation of the latter property is that for all $x\in X$ and all $b\in E$ there exists $a\in E$ such that $A_{ab}(x)=1$. 
\erem

\bdfn\label{d1_2014_12_15}
The incidence matrix $A$ and the random shift map $\hat\sg:\Om\to\Om$ alike, are said to be of finite range if for every $e\in E$ there exits a finite set $D_e\sbt E$ such that, for $m-a.e.\;  x\in X$,
$$
\{j\in E:A_{ej}(x)=1\}\sbt D_e
$$
Then also the random shift map $\hat\sg_A:\Om\to\Om$ is said to be of finite range. 
\edfn

\fr Note that finite range immediately entails local compactness of $m$ almost all symbol fibers $E_x^\infty$.
The following lemma is a reformulation of this condition.

\blem\label{l1} 
The random shift map $\hat\sg:\Om\to\Om$ is of finite range if and only if 
for every $l\in E$ there exists $\hat l \in E$ such that, for every $e>\hat l$,
$$
\sg_x^{-1} \left( [e]_{\th (x)}\right) \cap [0,...,l]_x =\emptyset 
\quad \text{for \ $m$-a.e.}\ x \in X. 
$$
\elem

\bdfn\label{d2_2014_12_15}
The incidence matrix $A$ and the random shift map $\hat\sg:\Om\to\Om$ alike, are said to be of bounded access if for all $b\in E$ there exists $b^*\in E$ such that for $m$-a.e. $x\in X$ there exists $a=a_x\le b^*$
such that $A_{ab}(x)=1$  (we recall that $E=\N$).
\edfn


\section{Preliminaries: Function Spaces and Bounded Distortion}
Various function spaces  will be used in the sequel. First we consider functions defined on a fiber $E_x^\infty$. The space of bounded 
continuous functions on $E_x^\infty$ will be denoted $\cC_b(E_x^\infty)$ and will be endowed with the supremum norm $\|\cdot\|_\infty$.
The $\al$--variation of a function $g:E_x^\infty \to \R$ is defined to be
$$
v_\al (g):=\inf \left\{ \frac{|g(\tau)-g(\om)|}{d(\tau,\om)^\al} \; : \;\; \tau\neq \om \; ,\;\; \tau_0=\om_0\right\}\,.
$$
The vector space $\cH_\al (E_x^\infty)$ by definition consists of all functions $g\in \cC_b (E_x^\infty)$ that have finite $\al$--variation. This space will be
endowed with the canonical norm 
$$
\|\cdot\|_\al = \|\cdot\|_\infty + v_\al(\cdot).
$$
Clearly, both normed spaces $\cC_b(E_x^\infty)$ and $\cH_\al (E_x^\infty)$ are Banach.

\

We now define global functions $g:\Omega \to \R$ and bring up their basic properties. 
Measurability of such functions  will always be with respect to 
the $\sg$--algebra $\cF\otimes|_\Om\cB$ on $\Om$.
A measurable function $g:\Om\to \R$
is called \emph{essentially continuous} if $g_x\in \cC_b (E_x^\infty)$ for $m$-a.e. $x\in X$, and if the essential supremum over $X$ of the measurable function $x\mapsto \|g_x\|_\infty $ is finite. The vector space of such functions will be denoted by $\cC_b(\Om)$. It becomes a Banach space when equipped with the norm
\beq\label{ess sup n}
|g|_\infty :=\ess\!\sup\{\| g_x\|_\infty:x\in X\}.
\eeq
We also have to consider various spaces of global H\"older functions especially in the context of the stochastic laws at the end of Section \ref{sec main thms}.
For the moment, let us introduce \emph{essentially $\al$-H\"older} functions. These will be measurable functions $g:\Omega \to \R$ such that $g_x\in \cH_\al (E_x^\infty )$ for $m$-a.e. $x\in X$ and such that 
$$
V_\al(g):= \ess\!\sup\{v_\a(g_x):x\in X\} <\infty \,.
$$
Functions that are in the same time essentially $\al$-H\"older and  essentially continuous form a Banach space,
denoted by $\cH_\al (\Om)$. The norm of a function $g\in \cH_\al (\Om)$ is
$$
|g|_{\al } := \ess\!\sup \{\|g_x\|_{\al}:x\in X\}\, .
$$

\

\fr For every integer $n\ge 1$ and every function $g:\Om\to\R$, let
$$
S_ng:=\sum_{j=0}^{n-1}g\circ \hat\sg^j,
$$
be the $n$th Birkhoff's sum of $g$ with respect to $\hat\sg$.
We will frequently need the following technical but indispensable fact.

\blem\label{2.9}
There exists a constant $C>0$ such that 
\beq\label{2.10}
\left| S_n g(x,\om) -S_ng(x,\tau)\right|\leq C V_\al (g) d(\om,\tau)^\al 
\eeq
and
\beq\label{2.11}
\left| \exp \left(S_n g(x,\om) -S_ng(x,\tau)\right)-1\right|\leq C V_\al (g) d(\om,\tau)^\al
\eeq
for every essentially $\al$-H\"older function $g:\Om\to\R$, $m$-a.e. every $x\in X$, every integer $n\ge 1$, and all $\om,\tau\in E_x^\infty$ with $\om|_{n-1}=\tau|_{n-1}$.
\elem

\fr The proof of this lemma is standard in hyperbolic dynamics and is left for the reader. 

\section{Transfer Operators and Fine Shifts}\label{Section_TO} 

We want to associate to a function $\ph:\Om\to \R$ a family of operators 
$\pf_x$, $x\in X$, by the formula
$$
\pf_xg_x(\om):=\sum_{e\in E\atop A_{e\om_0}(x)=1}g_x(e\om)e^{\ph_x(e\om)}\quad
\text{ where $g_x\in \cC_b\(E_x^\infty\)$ and $\om\in E_x^\infty$.} 
$$
Taking the particular function $\1_{[0,...,l]_x^c}$ where $l\in E$, then
$$\pf_x (\1_{[0,...,l]_x^c})(\om ) =\sum_{e >l\atop A_{e\om_0}(x)=1}e^{\ph_x(e\om)}
\quad \text{, $\om\in E_x^\infty$.}$$
Clearly, $\ph$ must satisfy some additional properties for these operators being well defined.

\bdfn\label{potential}
Given $\al>0$, we call an essentially $\al$-H\"older function $\ph:\Om\to\R$ a summable potential if
for every $e\in E$ there exists $0<c_e<C_e<\infty$ such that
\beq\label{ddd}
c_e\leq \exp {\ph _x} _{\big|[e]_x}\leq C_e
\eeq
and if the operators $\pf_x$ are uniformly summable in the following sense:
\beq \label{unif summ}
\lim_{l\to \infty}\left|  \pf_x \left( \1_{ [0,...,l]_x^c}\right)\right|_\infty =0\,.
\eeq
where the norm $|.|_\infty$ has been defined in \eqref{ess sup n}.
\edfn
\fr This summability condition deserves some
comments. First of all, \eqref{ddd} and \eqref{unif summ} imply that there exists $M<\infty$ such that 
\beq\label{const M}
\pf_x \1 \leq M \quad \text{for a.e. $x\in X$.}
\eeq
Then,  \eqref{unif summ} is formulated in terms of the transfer operator and thus this condition does depend on the
random incidence matrices. But, in general, one can use instead of  the above uniform summability  the following simpler condition which does only involve the function $\ph$ itself:
\beq\label{strongly summable}
\sum_{e\in E}\exp\(\ess\!\sup_{x\in X}{\ph _x} _{\big|[e]_x}\)<+\infty.
\eeq
Clearly, \eqref{strongly summable} implies  \eqref{unif summ}. However, 
 for some shifts and potentials \eqref{strongly summable} does not hold whereas the weaker
assumption \eqref{unif summ} does. This is the case for the following (deterministic) example.
Since we want also treat such shifts the condition \eqref{unif summ} is appropriate.

\bexam
The incidence matrix $A$ here does not depend on $x$ since we simple describe a deterministic example (which can be randomized in many ways). The one values of this matrix are defined as follows: $A_{1 0}=1$ and, for every $j\geq 1$, $A_{i j}=1$ if and only if 
$$i=(j-1) +n^{j+1} \quad \text{for some $n\in \N$.}$$
This shift has all required properties: it is mixing, of finite range and of bounded access (every deterministic shift is of bounded access). Consider then the potential
$$\ph (i)= -\ln (1+i) \quad , \quad i\in \N\,.$$
Clearly $\ph$ does not satisfy the summability condition \eqref{strongly summable} but it is obvious 
that  weaker summability condition \eqref{unif summ} does hold.
\eexam

\sp\fr Consider now the following three properties.

{\it
\sp\begin{itemize}
\item[(A)]
For every $e\in E$ there exists $M_e\in (0,+\infty)$ such that
$$
M_e^{-1}\le \pf_x\1|_{[e]_{\th(x)}}
\quad \text{for $m$-a.e. $x\in X$.}
$$
\item[(B)]
$\displaystyle
\lim_{e\to\infty}\ess\!\sup_{\!\!\!\!\! x\in X}\lt(\pf_x\1\big|_{[e]_{\th (x)}}\rt) =0.
$
\item[(C)] There are a number $0<\ka<1/4$ and a finite set $F\subset E$ such that
$$ 
\sup\lt(\pf_x \left( \1_{E_x^\infty \setminus [F]_x}\right)\rt)
\leq \kappa\inf\left(\pf_x \1|_{[F]_{\th (x)}}\rt)
\quad \text{for a.e. } x\in X\,.
$$
\end{itemize}
}

\sp\fr  These conditions deserve some comments. First of all,  as explained in the introduction, the goal here is to consider a situation where the transfer operator is no longer strictly positive in the sense that $ \inf \pfx \1>0$.
Clearly, the condition (B) is responsible for this whereas (C) gives some control of how 
$\pf_x\1 (y)\to 0$ as $y_0\to\infty$ and (A) is a weak lower bound.
An other remark is that these conditions rely only on the shift and the potential itself and so no higher iterates are involved. They are, at least for reasonable potentials, quite simple to check.
Indeed, condition (C) relies on the potential but, for summable potentials, (A) and (B) totally rely on the shift, hence the incidence matrix. 

\blem Let $\ph:\Omega \to \R$ be a summable potential. Then we have the following:
\ben
\item The shift is of bounded access if an only if (A) holds.
\item The shift is of finite range if and only if (B) holds.
\een
\elem

\fr The proof of this
remark is obvious.

\bdfn\label{fine}
A random shift $\hat\sg:\Om\to\Om$ along with a potential $\ph : \Om\to \Om$
is called fine if $\ph$ is summable and if $(A,B,C)$ hold or, equivalently, if $\ph$ is summable,
the shift $\hat\sg$ is of finite range and of bounded access and if (C) holds.
\edfn

\sp Here are some more properties that will be useful in the sequel.

The condition \eqref{ddd} is just a weak bound and it implies the following immediate observation.

\blem\label{l2}
Suppose that the random shift map $\hat\sg:\Om\to\Om$ is of finite range and that $\ph:\Om\to\R$ is essentially 
$\a$-H\"older and satisfies \eqref{ddd}. Then, for every $e\in E$ and every integer $n\geq 1$ there exists $c^*=c^*(e,n)>0$ such that 
$$
\exp\left(\inf\(S_n \ph|_{\{x\}\times [e]_x}\)\right)\geq c^* \quad \text{for $m$ a.e.} \  x \in X.
$$
\elem

\

\fr For a summable potential it is clear that $\pf_xg_x\in \cC_b\(E_{\th(x)}^\infty\)$ and so 
$$
\pf_x\(\cC_b\(E_x^\infty\)\)\sbt \cC_b\(E_{\th(x)}^\infty\).
$$
We record the following.

\blem\label{2.0.1}
If $\ph:\Om\to\R$ is a summable potential, then for each $x\in X$, the linear operator 
$$
\pf_x:\cC_b\(E_x^\infty\)\to \cC_b\(E_{\th(x)}^\infty\)
$$ 
is bounded and its norm is bounded above by $M$. Also $\pf_x\(\cH_\al\(E_x^\infty\)\)\sbt \cH_\al (E_{\th(x)}^\infty\)$ and the linear operator 
$
\pf_x:\cH_\al\(E_x^\infty\)\to \cH_\al\(E_{\th(x)}^\infty\)
$ 
is bounded.
\elem

\bpf
The first assertion of this lemma is obvious while the second results by a standard calculation.
\epf

\fr The operators $\pf_x$, $x\in X$, are frequently referred to as fiberwise transfer operators. They give rise to the global operators defined as follows. If $g\in \cC_b(\Om)$ and $x\in X$, we set
$$
(\pf g)_x:=\pf_{\th^{-1}(x)}g_{\th^{-1}(x)}\in \cC_b\(E_x^\infty\).
$$

\blem\label{o1_2015_01_16}
If $\ph:\Om\to\R$ is a summable potential, then the function $X\ni x\mapsto \pf_x$ is measurable in the sense that for each $g\in \cC_b(\Om)$ the function
$$
\Om\ni\(x,\om)\mapsto \pf_xg_x(\om)\in\R
$$
is measurable. In consequence $\pf g\in\cC_b(\Om)$.
\elem

\bpf
Note that, given $e\in E$, both functions $\Om\ni(x,\om)\mapsto g_x(e\om)$ and $\Om\ni(x,\om)\mapsto e^{\ph_x(e\om)}$ are measurable and can be extended to measurable functions on $X\times E^\N$ by putting the value zero outside $\Om$. Since we can extend measurably in the same way the incidence matrix $A$,
the function
$$
X\times E^\N \ni(x,\om)\mapsto A_{e\om_0}(x)  g_x(e\om)e^{\ph_x(e\om)}
$$
is measurable too. This shows measurability of the function $X\ni x\mapsto \pf_x(g)$ since the former one can now be expressed as a convergent series of measurable functions.
\epf

\fr As a direct consequence of this observation and of Lemma~\ref{2.0.1}, we get the following. 

\blem\label{2.0.1b}
If $\ph:\Om\to\R$ is a summable potential, then the linear operator 
$$
\pf:\cC_b(\Om)\to\cC_b(\Om)
$$
is bounded and its norm is bounded above by $M$. Also $\pf\(\cH_\al (\Om))\sbt \cH_\al (\Om)$ and the linear operator 
$
\pf:\cH_\al (\Om)\to \cH_\al (\Om)
$
is bounded.
\elem

\sp\fr For $n\geq 1$, we define the iterated operator
$$
\pf_x^n := \pf_{\th^{n-1}(x)}\circ...\circ \pf_x:\cC_b\(E_x^\infty\)\to \cC_b\(E_{\th^n(x)}^\infty\)
$$
Of course
$
\pf_x^n \(\cH_\al\(E_x^\infty\)\)\sbt \cH_\al (E_{\th^n(x)}^\infty\).
$
A standard straightforward inductive calculation shows that
$$
\pf_x^n(g)(\om)=\sum_{\tau\in E_x^n\atop A_{\tau_{n-1}\om_o}(x)=1} \exp\(S_n\ph(x,\tau\om)\)g(\tau\om).
$$
As an immediate consequence of Lemma~\ref{2.9} we get the following.

\blem\label{2.12}
For every $n\geq 1$, every $x\in X$, and all $\om,\tau\in E_x^\infty$, with $\om_0=\tau_0$,
$$
\frac{\pf^n_x\1 (\om)}{\pf^n_x\1 (\tau)}
\leq 1+ CV_\al(\ph)d(\om,\tau)^\al\, .
$$
In particular, 
$$
\pf_x^n\1 (\om)\leq K \pf_x^n\1 (\tau),
$$
where $K= 1+ CV_\al (\ph )$.
\elem

Our goal now is three-folded: to prove the existence of conformal measures, 
of their invariant versions and 
finally to obtain exponential rate of convergence for the iterated normalized operator along with stochastic laws. 
This will be done for random shifts and potentials that satisfy the conditions formulated in the next section.


\section{Random Measures and Main Theorems}\label{sec main thms}
 
The first thing we want to do in this section is to recall the concept of random measures and to bring up some of its basic properties. We do this in our context of the measurable space $\Om$ and measure $m$ on $X$. 

\bdfn\label{random measures}
A measure $\nu $ on $\(\Om, \cF\otimes|_\Om\cB\)$ with marginal $m$, i.e. such that  $$
\nu\circ\pi_X^{-1}= m,
$$
is called  a random measure if its disintegrations $\nu_x$, $x\in X$, respectively belong to the spaces $\cP(E_x^\infty)$ of probability measures on $\(E_x^\infty,\cB_x\)$. The space of all random measures (on $\Om$ with marginal $m$) will be denoted by $\cP_m(\Om)$. Of course every element of $\cP_m(\Om)$ is a probability measure on $\(\Om, \cF\otimes|_\Om\cB\)$
\edfn

\fr According to this definition, for every random measure $\nu$, i. e. belonging to $\cP_m(\Om)$, and for every every function $g\in \cC_b(\Om)$, we have
\beq\label{1_2015_01_17}
\nu(g)
 :=\int_\Om g\,d\nu
 =\int_X \int _{E_x^\infty} g_x\, d\nu_x\, dm(x)
\leq |g|_\infty.
\eeq
This formula naturally introduces the space $g\in L^1(\nu)$, i.e. the space of all real-valued measurable functions on $\Om$ integrable with respect to $\nu$. Precisely, $g\in L^1(\nu)$ if and only if 
$$
\int_X \int _{E_x^\infty} |g_x|\, d\nu_x\, dm(x) <+\infty\,.
$$

\sp We would like to mention that random measures, as defined in Crauel's book \cite{Cra02}, are all probability measures on the set $X\times E^\N$ with marginal $m$ and distintegration measures being probabilities in $E^\N$. Denote them just by $\cP_m$. But in this paper we actually are interested only in the class $\cP_m(\Om)$, defined few paragraphs above, i.e. in the measures in $\cP_m$ whose supports lie $\Om$, i.e. such that $\nu(\Om)=1$. 

\sp The key concept pertaining to random measures is that of narrow topology (a version of weak convergence). Namely, if $\La$ is a directed set, then a net $\(\nu^\a\)_{\a\in\La}$ in $\cP_m$ is said to converge to a random measure $\nu\in \cP_m$ if
$$
\lim_{\a\in\La}\nu^{\a} (g) =\nu (g) \quad \text{for every $g\in \cC_b(X\times E^\N)$. }
$$ 
This concept of convergence defines a topology on $\cP_m$ called in \cite{Cra02} the narrow topology. The narrow topology on $\cP_m(\Om)$ is the one inherited from the narrow topology on $\cP_m$. Since $\1_\Om$, the characteristic function of $\Om$ in $X\times E^\N$, belongs to $\cC_b(X\times E^\N)$, we have that
$$
\nu(\1_\Om)=\lim_{\a\in\La}\nu^{\a} (\1_\Om) =1
$$ 
for any net $\(\nu^\a\)_{\a\in\La}$ in $\cP_m(\Om)$ converging to a random measure $\nu\in \cP_m$. This means that then $\nu\in \cP_m(\Om)$, leading to the following.

\bprop\label{p2_2015_01_17}
$\cP_m(\Om)$ is a closed subset of $\cP_m$ with respect to the narrow topology on $\cP_m$.
\eprop

\fr Recall from \cite{Cra02} that a subset $\Ga$ of $\cP_m$ is called tight if its projection $\pi_{E^\N}^{-1}(\Ga)$ on $E^\N$ is a tight subset of Borel probability measures on  $E^\N$, the latter (commonly) meaning that for every $\ep>0$ there exists a compact set $K_\ep\sbt E^\N$ such that $\nu\circ\pi_{E^\N}^{-1}(K_\ep)\ge 1-\ep$ for all $\nu\in\Ga$. For us, the crucial property of narrow topology is that of Prohorov's Compactness Theorem (Theorem 4.4 in \cite{Cra02}) which asserts that a subset $\cM \subset \cP_m$ is relatively compact with respect to the narrow topology if and only if it is tight. Along with Proposition~\ref{p2_2015_01_17}, this entails the following.

\bthm\label{t3_2015_01_17}
A subset $\Ga\subset \cP_m(\Om)$ is relatively compact with respect to the narrow topology if and only if it is tight. Furthermore, it is compact if and only if it is tight and closed.
\ethm

\fr We will use tightness heavily in the next section and for this we will need Proposition~4.3 from \cite{Cra02}), which provides convenient sufficient conditions for tightness to hold. For this in turn, we will need concepts of random closed and compact sets. Indeed, following Definition~2.1 in \cite{Cra02}) we say that a function
$$
X\ni x\mapsto C_x,
$$
ascribing to each point $x\in X$ a closed subset $C_x$ of $E^\N$, is called a random closed set if for each $\om\in E^\N$ the function
$$
X\ni x\mapsto \dist(\om,C_x)\in\R
$$
is measurable. Since the probability measure $m$ on $X$ is assumed to be complete, being a closed random set precisely means (see Proposition~2.4 in \cite{Cra02}) that the union
$$
\bu_{x\in X}\{x\}\times C_x \  \text{ (all sets $C_x$ are assumed to be closed)}
$$
is a measurable subset of $X\times E^\N$. A random closed set $X\ni x\mapsto C_x$ is called a random compact set if all sets $C_x$, $x\in X$, are compact (in $E^\N$). A function $X\ni x\mapsto V_x$ is called a random open set if the function $X\ni x\mapsto E^\N\sms V_x$ is a
closed random set. The theorems about tightness announced above are these (see Proposition~4.3 in \cite{Cra02}).

\bprop\label{p4_2015_01_17}
A subset $\Ga\subset \cP_m(\Om)$ is tight if and only if for every $\ep>0$ there exists a random compact set $X\ni x\mapsto K_x$ such that $K_x\sbt E_x^\infty$ for all $x\in X$ and
$$
\int_X\nu_x(K_x)\, dm(x)\ge 1-\ep
$$
for all $\nu\in\Ga$.
\eprop
\fr As an immediate consequence of this proposition, we get the following. 

\bcor\label{c5_2015_01_17}
Let $\Ga$ be a subset of $\cP_m(\Om)$. Suppose that for every $\ep>0$ there exists a random compact set $X\ni x\mapsto K_x$ such that $K_x\sbt E_x^\infty$ for all $x\in X$ and
$$
\nu_x(K_x)\ge 1-\ep
$$
for all $\nu\in\Ga$ and for $m$-a.e. $x\in X$, then $\Gamma$ is tight.
\ecor

Passing to dynamics, i. e. to the fine random shift $\hat\sg:\Om\to\Om$ and fine potential $\ph:\Om\to\R$, as indicated in the introduction, the measures we are looking for are defined in dynamical terms and form some special random measures on $\Om$.

\bdfn\label{d120141216}
A random conformal measure is a measure $\nu=(\nu_x)_{x\in X} \in \cP_m(\Om)$ for which there exists a measurable function $X\ni x\mapsto \l_x\in (0,+\infty)$ such that
\beq \label{eq conf m}
\pf_x^* \nu_{\th (x) } = \l_x \nu_x \quad \text{for $m$-a.e.} \ x \in X\,.
\eeq
\edfn


\fr We can now present the main results of this paper.

\bthm \label{thm 1}
Let the random shift $\hat\sg:\Om\to\Om$ and the potential $\ph:\Om\to\R$ be fine.
Then there exists a random conformal measure, i.e. a measure $(\nu_x)_{x\in X}$ that satisfies
\eqref{eq conf m} for some measurable function  $X\ni x\mapsto \l_x\in (0,+\infty)$.
In addition, 
$$
\| \log \l \|_\infty <\infty\, .
$$
\ethm


\fr Let
$$
\npf_x:=\l_x^{-1}\pf_x:\cC_b\(E_x^\infty\)\to \cC_b\(E_{\th(x)}^\infty\)
$$
and let
$$
\npf_x^n := \npf_{\th^{n-1}(x)}\circ...\circ \npf_x:\cC_b\(E_x^\infty\)\to \cC_b\(E_{\th^n(x)}^\infty\)
$$
Our second main theorem is this.

\bthm \label{thm 2}
Let the random shift $\hat\sg:\Om\to\Om$ and the potential $\ph:\Om\to\R$ be fine.
Then
\ben
\item There exists a unique positive function $\rho\in \cH_\al(\Om)$ 
such that $\npf\rho=\rho$, which fiberwise means that 
$$
\npf_x\rho_x =\rho_{\th (x)}\quad \text{for $m$-a.e. $x\in X$.}
$$ 
\item There exist constants $B>0$ and $\vth\in (0,1)$ such that for $m$-a.e. $x\in X$
$$ 
\left\| \npf_x^n g - \nu_x(g) \rho _{\th^n (x)} \right\|_\al \leq B\vth ^n \quad \text{for every} \ g\in \cH_\al (E_x^\infty) \,.
$$
\een
\ethm

\brem
The Condition (C) is only required for the existence of conformal measures. Consequently, if for some reasons there already exists a conformal measures (see, for example, the discussion on interval maps in Section 6.2), then Theorem \ref{thm 2} does hold without the assumption that Condition (C) holds.
\erem

\brem\label{quenched}
We would like to remark that both Theorems\, \ref{thm 1} and \ref{thm 2} belong to quenched type of results; see \cite{ANV} for an extensive discussion of this concept.
Compare also Remark~\ref{annealed} on annealed results.
\erem

\sp\fr As an immediate consequence of part (1) of this theorem, we get the following.

\bcor\label{c6_2015_01_17}
Let the random shift $\hat\sg:\Om\to\Om$ and the potential $\ph:\Om\to\R$ be fine.
Let $\rho\in \cH_\al(\Om)$ come from Theorem~\ref{thm 2}. Then for the fiber measures $\mu_x=\rho_x\nu_x\in \cP(E_x^\infty)$, $x\in X$, we have that
$$
\mu_x\circ\sg_x^{-1}=\mu_{\th(x)},
$$
and $\mu$, the random measure on $\Om$ with disintegration $(\mu_x)_{x\in X}$, is the unique $\hat\sg$-invariant measure on $\Om$ absolutely continuous with respect to $\nu$. In addition, the two measures $\mu$ and $\nu$ are equivalent.
\ecor

We want to conclude this section with two striking stochastic consequences of Theorem~\ref{thm 1} and  Theorem~\ref{thm 2}, particularly its item (2). We mean exponential decay of correlations and the Central Limit Theorem. These follow from Theorems~\ref{thm 1} and \ref{thm 2} in an analogous way as corresponding theorems in \cite{MyUrb2014} followed from its respective counterparts of Theorems~\ref{thm 1} and \ref{thm 2}. We therefore only formulate the exponential decay of correlations and the Central Limit Theorem inviting the interested readers to look for proofs to \cite{MyUrb2014}. In order to formulate these theorems we need some new, more general, function spaces.

\sp Given a number $0<p\le \infty$ let $\cH_\a^p (\Om)$ be the space of functions $g:\Om\to \R$ with H\"older fibers $g_x\in \cH_\a (E_x^\infty)$ and such that
$\|g_x\|_\a\in L^p(m)$ for $m$-a.e. $x\in X$. The canonical norm on this spaces is 
$$
|g|_{\a , p}= \left( \int_X  \|g_x\|_\a ^p\, dm(x)\right)^\frac{1}{p}\,.
$$
It makes $\cH_\a^p (\Om)$ a Banach space.
Replacing in this definition the $\a$--H\"older condition on the fiber $E_x^\infty$ by a $L^1 (\nu_x)$ condition leads to a space of functions that will be denoted by $L_{\nu }^{ 1, p}(\Om)$. The natural norm is in this case 
$$
| g|_{\nu }^{ 1, p} = \left(\int _X \| g_x\|^p _{L^{1} (\nu_x )}\, dm(x)\right)^\frac1p \,.$$
Clearly, if $p=1$ then $L_{\nu }^{ 1, 1}(\Om)=L^1 (\nu)$. 
In both cases we also consider $p=\infty$ and then the $L^p$ norms are understood as the sup--norms. 

\bthm[Exponential Decay of Correlations]\label{thm-correlations}
Let the random shift $\hat\sg:\Om\to\Om$ and the potential $\ph:\Om\to\R$ be fine.
Let $\mu$ be the corresponding $\hat\sg$-invariant measure on $\Om$ produced in Theorem~\ref{thm 2}. Let $p,q\in [1,\infty ]$ be such that $\frac1p+\frac1q=1$.
Then, for every  $g\in L_{\nu }^{ 1, p}(\Om)$, $h\in \cH_\a ^q (\Om)$ with $\int_{X_x^\infty}h_x\,d\mu_x=0$ and for every $n\geq 1$, we have
\begin{align*}
\left| \int_{\Om}(g\circ\hat\sg^n) \, h \, d\mu \right| 
= \left| \int _X \int _{E_x^\infty}(g_{\shift^n(x) }\circ\hat\sg_x^n) \, h_x\, d\mu_x\, dm(x) \right|
\leq b \vartheta ^n \; |g|_{\nu}^{ 1, p} \;  |h|_{\b,q}
\end{align*}
for some positive constant $b$ and some $\vartheta\in(0,1)$. 
\ethm

\bthm[Central Limit Theorem]\label{CLT}
Let the random shift $\hat\sg:\Om\to\Om$ and the potential $\ph:\Om\to\R$ be fine.
Let $\mu$ be the corresponding $\hat\sg$-invariant measure on $\Om$ produced in Theorem~\ref{thm 2}. Let $\psi \in \cH_\a (\Om)$ such that
$ \int _{E_x^\infty} \psi _x d\mu_x =0$, $x\in X\,.$
If $\psi$ is not cohomologous to $0$, then there exists $\sg>0$ such that, for every $t\in \R$, 
$$
\mu \left( \lt\{z\in \Om: \; \frac1{\sqrt{n}} S_n \psi (z) \leq t\rt\}\right) 
\longrightarrow
\frac1{\sg \sqrt{2\pi }}\int _{-\infty}^t exp (-u^2/2\sg^2 ) \, du.
$$
\ethm

\

\brem \label{annealed}
We would like to remark that both Theorems\, \ref{thm-correlations} and \ref{CLT} belong to annealed type of results; see \cite{ANV} for an extensive discussion of this concept.
Compare also Remark~\ref{quenched} on quenched results.
\erem

\section{Examples}\label{sec examples}

 The present work is related to  our former work on transcendental dynamics \cite{MyUrb2014}. In contrast to the case of hyperbolic rational functions, the dynamics of a general hyperbolic transcendental function on the Julia set can not be represented by a symbol dynamics. 
However,  some particular entire functions act on a dynamically important subset of the Julia set, the so called set of "landing-" or "end-"points, in a similar way as countable Markov shifts that we consider in here
(see \cite{Baran07, BaranKarp07}).

\subsection{Non-homogenous random walks.} There is a natural relation between random dynamics and stochastic processes in random environments (see \cite{Kif98, ANV}) and the particular case of a nearest neighbor random walk totally fits into our setting. Indeed, consider the random walk on $\Z$ given by 
$X_n=\sum_{j=1}^ n Y_j$ where $Y_j$ are iid random variables having equiprobably chosen values in $\{-\eta, ...,-1,0,1, ... , \eta\}$,
$\eta\geq 1$ being the number of authorized neighbors. Such a walk is equivalent to the shift map whose incidence matrix has coefficient $A_{i j}=1$ if and only if $|i-j|\leq \eta$. Note that such shifts satisfy all the required properties in our present paper: these are mixing, of bounded access, and of finite range. 

Our random setting can also represent random walks where the equiprobable choice described above is replaced by transition laws that depend on the site and also on time. Such random walks are called inhomogenous random walks in random environments \cite{Zeitouni04}

Related to these examples are the random walks with random potentials considered in \cite{RSY2013}. For appropriate choices of the potential $V$ in \cite{RSY2013} we are again in the scope of the present paper. To the following more general example all the results of the present paper do apply: consider a random walk in $\Z$ where the number of neighbors of a site $i$ is some function $\eta (i)\in \N\setminus \{0\}$. When the walk is on site $i$ then, with probability $p_i=\frac1{2 \eta(i)+1}$, one of the sites in 
$\{i-\eta (i),...,0,1,..., i+\eta (i)\}$ is chosen. Such a walk corresponds to a deterministic shift with incidence matrix determined by $A_{i j}=1$ if and only if $|i-j|\leq \eta (i)$. The natural potential is then
$$\ph (\om ) = \log p_{\om _0}$$
or even $\ph_\b  (\om ) = \b \log p_{\om _0}$ for some $\b \in \R$.
For simplicity let us just consider the particular choice 
$$\b =1 \quad \text{and} \quad \eta (i)= 4^{|i|}\,.$$
Now, a simple calculation shows that this potential satisfies all the required properties of the present paper.
In particular, the condition (C) is satisfied with $F=\{-1,0,1\}$ and $\kappa =1/5$. This shift is a fine shift
and can be turnt into a random shift in various ways. An interesting random shift emerges when the function $\eta$
is replaced by a random variable $\eta _x(i)$.

\subsection{Interval maps and iterated function systems.}  
Another field of application of our symbol considerations are interval maps and iterated function systems, or even graph directed Markov systems. Let us just indicate that various examples of the later are contained in 
\cite{MauldinUrbanski} (deterministic iterated function systems) and in \cite{RoyUrb2011} (random graph directed Markov systems). Concerning interval maps, Example 3.3 in \cite{GI05} is a typical 
example whose associated shift is of bounded access and finite range (although it is again deterministic but it can obviously be perturbed to a random map). Note that here, like in many other cases of interval maps, there exists for free a conformal one, Lebesgue measure. The reason is that the associated Markov partition is a partition of the whole unit interval. Since in the present paper we construct conformal measures we can consider interval maps whose underlying partition does not cover the whole unit interval.


\section{The Existence of Conformal Measures}

This section is devoted to prove Theorem~\ref{thm 1}. It follows from the invariance relation \eqref{eq conf m} that 
$$
\l_x = \int \cL_x \1 \, d\nu_{\shift (x)},
$$ 
and so our task is to look for random measures
$(\nu_x)_{x\in X}$ that are invariant under the map $\Phi:\cP_m(\Om) \to\cP_m(\Om)$ that is fiberwise defined as follows:
$$
\Phi(\nu)_x:=\frac{\pf_{x}^*\nu_{\th(x)}}{\pf_{x}^*\nu_{\th(x)}(\1)}.
$$
We want to obtain these measures in the usual way by employing Schauder--Tychonov Fixed Point Theorem.
But, since the sets $E_x^\infty$ need not be compact, this can be done only if a convex compact and $\Phi$--invariant 
subset $\cM$ of $\cP_m(\Om)$ can found, and if in addition the map $\Phi$, restricted to this set, is continuous. 
Such a subspace will now be found with, in particular, the help of the assumption (C). Remember that this assumption involves a finite union of cylinders $F$. Re-enumerating $E=\N$ if necessary, we may assume that $F=\{0,1,...,q\}$ for some $q\geq 0$. We shall prove the following.

\blem\label{3.2} 
Let the random shift $\hat\sg:\Om\to\Om$ and the potential $\ph:\Om\to\R$ be fine.
Then, there exists $(\ga_n)_{n=0}^\infty$, a sequence of measurable positive functions defined on $X$, with the following properties:

\sp\ben
\item $\ga_0 (x) \equiv \frac12$,

\sp\item $\lim_{n\to \infty } \ga_n (x) =0$\, for $m$-a.e. $x\in X$,

\sp\item
there exists an ascending sequence of 
compact random sets, whose fibers we denote by $K_{n,x}$, $x\in X$, 
and a measurable point $x\mapsto \xi (x)$
such that 

\sp\ben
\item[(3.1)] $\xi (x) \in K_{0,x}\subset [F]_x$ for $m$-a.e. $x\in X$,
and   

\sp\item[(3.2)]the set 
$$
\cM :=\big\{  \nu \in \cP_m(\Om) \; , \;\; \nu_x(K_{n,x}^c ) \leq \ga_n (x) \;\; \text{for all}\;\;
n\geq 0 \text{ and a.e. } x\in X\big\}
$$
is $\Phi$--invariant, precisely meaning that $\Phi(\cM) \subset \cM$.
\een
\een
\elem

\bpf First of all, notice that $\displaystyle \F= \bigcup_{x\in X}\{x\}\times [F]_x$ is a closed random set
with  non-empty fibers $[F]_x$ (see Remark \ref{2.001}). Therefore, the Selection Theorem (Theorem 2.6
in \cite{Cra02}) implies that there exists a measurable point $\xi\in \F$, i.e. a measurable map
$x\mapsto \xi (x) \in [F]_x$.

We will now define inductively a required sequence $(\ga_n)_{n=0}^\infty$. During this inductive process we will also define an auxiliary sequence $\(X\ni x\mapsto N_n(x)\in\N\)_{n=0}^\infty$ of integer-valued functions meeting the following properties:

\ben

\sp\item the functions $N_n(x)$ are measurable for all integers $n\geq 0$.

\sp\item $\(N_n(x)\)_{n=0}^\infty$ is an increasing sequence for every $x\in X$.

\sp\item $N_0\equiv q$ and $N_n(x)\geq \xi_n (x)$ for all integers $n\geq 0$ and all $x\in X$,
\een
where $\xi_n (x)$ is the $n$--th coordinate of the measurable function $\xi$ define just above. 

Assuming that such a sequence $ N_n(x)$ is given, define further
$$
N_{j,n-j} (x):= N_n (x) 
$$ 
for all $j=0,1,...,n$, and define also
\beq \label{3.1}
K_{n,x}:= \big\{\om\in E_x^\infty: \,\om_i \leq N_{i,n} (\th^i (x)) \; \; \text{for all}\;\; i\geq 0 \big\}.
\eeq
Notice that these sets $K_{n,\om}$ are compact subsets of $E_x^\infty$ and that, due to the measurability of the functions $N_n$, the family $(K_{n,x})_{x\in X}$ forms a compact random set. Moreover,  
$$
\xi (x) \in K_{0,x}\subset [F]_x\quad\text{and}\quad K_{n,x} \subset K_{n+1,x}
$$
for all integers $n\geq 0$ and all $x\in X$.
In order to demonstrate $\Phi$--invariance of the set $\cM$, let us first make the following observation.
Let $\tau\in K_{n+1,\th (x)}$ and $j\in E$ be such that $j\tau\in E_x^\infty$ or, equivalently, such that $A_{j\tau_0}(x)=1$.
Then $ \1_{K_{n,x}^c}(j\tau)= 1$ if and only if either
$$
j> N_{n} (x) \quad \text{or else}\quad \tau_i> N_{i+1,n}(\th^{i+1} (x) ) = N_{n+i+1}(\th^{i+1} (x) ) \; \text{ for some $i \geq 0$}\,.
$$
On the other hand, 
$\tau_i\leq N_{i,n+1}(\th ^i ( \th (x))) = N_{n+i+1}(\th ^{i+1} (x) )$ since $\tau\in K_{n+1,\th (x)}$. 
Consequently, $\tau\in K_{n+1,\th (x)}$ and $ \1_{K_{n,x}^c} (j\tau)= 1$ if and only if  $j> N_{n} (x)$. This allows us to make the following estimation valid for all $\tau\in K_{n+1,\th (x)}$:
\beq\label{1_2014_12_16}
\pf_x\1_{K_{n,x}^c }(\tau)
=\sum_{j\in E\atop j\tau\in E_x^\infty} e^{\ph_x (j\tau) } \1_{K_{n,x}^c} (j\tau)
=\sum_{j> N_{n} (x) \atop j\tau\in E_x^\infty} e^{\ph _x (j\tau)} 
\leq Z_x (N_n (x))
\eeq
where, for every $k\geq 0$ and every  $x\in X$,
$$
Z_x(k):= \sup \lt(\pf_x\1_{[0,1,...,k]_v^c}\rt) 
=\sup _{\tau\in E_{\th (x)}}\lt(
\sum_{j> k \atop j\tau\in E_x^\infty} e^{\ph _x (j\tau)}\rt)\,.
$$
Summability of the potential $\ph$ (see Definition~\ref{potential}) implies that
\beq \label{3.7} 
\lim_{k\to \infty} Z_x(k) =0 \  \text{ for } \ m-a.e.\  x\in X.
\eeq
Take now an arbitrary $\nu\in \cM$. Since
\begin{align*}
\pf_x^* \nu_{\th (x)} (K_{n,\om }^c ) 
& = \int \pf_x\1_{K_{n,x}^c }d\nu_{\th (x)}\\
&=\int _{K_{n+1,\th (x) }^c}\pf_x \1_{K_{n,x}^c }d\nu_{\th (x)}+\int_{K_{n+1,\th (x)}} \pf_x \1_{K_{n,x}^c }d\nu_{\th (x)}\, ,
\end{align*}
it follows from summability of $\ph$  and from \eqref{1_2014_12_16} that
\beq\label{3.5}
\begin{aligned}
\pf_x^* \nu_{\th (x)} (K_{n,x}^c )
&\leq\| \pf_x \1\|_\infty \nu_{\th (x)}(K_{n+1,\th (x) }^c)
+ Z_x (N_n (x)) \nu_{\th (x)}(K_{n+1,\th (x)})\,\\
&\leq  M \ga _{n+1}(\th (x)) + Z_x (N_n (x)).
\end{aligned}
\eeq
On the other hand, since $\nu\in \cM$, we have
$$
\begin{aligned}
\pf_x^*\nu_{\th (x)} (\1 ) 
&\geq \int _{K_{0 , \th (x) }}\pf_x\1 d \nu _{\th (x)}
\geq \nu_{\th ( x )} (K_{0 , \th ( x ) }) \inf_{K_{0 , \th ( x ) }}\(\pf_x \1\)\\ 
&\geq \frac12\inf _{K_{0 , \th (x) }}\(\pf_x\1\) \\
&\geq \frac12 \inf _{[F]_{\th (x)}}\(\pf_x\1\):=c>0
\end{aligned}
$$
For some constant $c>0$ because of (A).
Together with \eqref{3.5} this leads to
$$
\Phi (\nu )_x \left(K_{n,x}^c \right)\leq 
\frac1{c} \Big(
M \ga _{n+1}(\th (x)) + Z_x (N_n (x))\Big)\,.
$$
Therefore, assuming that a measurable function $\g_n$, is given, if we can find measurable functions $\ga_{n+1}$
 and $ N_n$ such that 
\beq\label{3.6}
M\ga _{n+1}(\th (x)) + Z_x(N_n (x))\leq c \ga _n (x), 
\eeq
then the proof of our lemma will be complete. 

\sp Let us first consider the case of $n=0$. Since we have put $\g_0=1/2$ and $N_0=q$,  formula \eqref{3.6} then takes on the form
$$
M\g_1(\th (x))+Z_x(q)\le \frac14\inf_{[F]_{\th (x)}}\(\pf_x\1\)\,.
$$
Therefore, condition (C) yields that it is sufficient to take measurable
$\g_1(\th (x))$ such that 
$$
M\g_1(\th (x)) < \left(\frac14-\kappa\right)\inf _{[F]_{\th (x)}}\(\pf_x\1)\,.
$$
This is possible as $\kappa<1/4$.

\sp Suppose finally that, for some $n\ge 0$, a measurable function $\g_n$, is given. Then, because of \eqref{3.7}, one can find a measurable function $N_n\geq \xi_n$ such that $Z_x(N_n (x))\leq \frac12 c \ga_n (x)$.
So, setting then
$$
\g_{n+1}(x):=M^{-1}\(c\ga_n(\th^{-1}(x))-Z_{\th^{-1}(x)}(N_n (\th^{-1}(x))\)
$$
defines a measurable function from $X$ to $(0,+\infty)$ for which \eqref{3.6} holds. The proof is complete.
\epf


\fr The set of measures $\cM$ produced in this lemma is not only $\Phi$-invariant but it has all the required properties.

\blem \label{c.2}
The set $\cM$ is non-empty, convex, and compact.
\elem

\bpf
With the same notation as in Lemma~\ref{3.2}, we have
$\xi (x) \in K_{0,x}\subset K_{n,x}$ for every $n\geq 0$ and a.e. $x\in X$.
Consider then the measure $\nu$ defined fiberwise by $\nu_x= \delta_{\xi (x)}$,
the Dirac $\delta$--measure supported at 
$\xi (x)$, $x\in X$. Then obviously $\nu_x(K_{n,x}^c )=0\leq \ga_n (x)$,
which shows that $\nu \in \cM$ and thus $\cM\neq \emptyset$. Convexity of $\cM$ is obvious.
In order to prove compactness of $\cM$, we shall show that $\cM$ is closed and tight. 

Tightness first. Fix $\ep >0$ arbitrary.
Since, by item (2) of Lemma~\ref{3.2}, the sequence $(\ga_n)_n$ converges pointwise to $0$, there exists an integer $k=k(\ep)\ge 0$ such that $m\big(\ga_k ^{-1} ([\ep/2 ,+\infty ))\big)<\ep/2$.
Then, for every $\nu\in \cM$ we have that 
\begin{align*} 
\int _X \nu_\om \left( K_{k,x}^c\right) dm(x) 
&=\int _{\ga_k ^{-1} ([\ep/2 , \infty ))} \!\!\! \nu_\om \left( K_{k,x}^c\right) dm(x)
   +\int _{\ga_k ^{-1} ([0,\ep/2 ))} \!\!\! \nu_\om \left( K_{k,x}^c\right) dm(x)\\
&\leq m(\ga_k ^{-1} ([\ep/2 , \infty ))) + \frac{\ep}2  \\
&<\frac{\ep}2+\frac{\ep}2 
=\ep \,.
\end{align*}
Therefore, the set $\cM$ is tight (see Proposition~\ref{p4_2015_01_17}). 

To see that $\cM$ is closed let $\Lambda$ be a directed set and $(\nu^\a)_{\a\in\La}$ a net in $\cM$ converging to a measure  $\nu \in  \cP_m(\Om)$ in the narrow topology. Let $H$ be an arbitrary measurable subset of $X$. Then, for every integer $n\ge 0$, the function
$$
X\ni x\mapsto 
\begin{cases}
K_{n,x}^c   &\text{ if }  \  x\in H \\
\es  &\text{ if }  \  x\notin H 
\end{cases}
$$
defines a random open set. It then follows from Portmenteau's Theorem (see Theorem~3.14 (iv) in \cite{Cra02}) that
$$
\begin{aligned}
\int_H\nu_x\(K_{n,x}^c\)\,dm(x)
&=\nu\lt(\bu_{x\in H}\{x\}\times K_{n,x}^c \rt)
\le \varliminf_{\a\in\La}\nu_\a\lt(\bu_{x\in H}\{x\}\times K_{n,x}^c \rt) \\
&=\varliminf_{\a\in\La}\int_H\nu_{\a,x}\(K_{n,x}^c\)\, dm(x) \\
&\le \int_H\g_n(x)\, dm(x).
\end{aligned}
$$
Hence, arbitrariness of $H$ yields $\nu_x\(K_{n,x}^c\)\le \g_n(x)$ for $m$-a.e. $x\in X$. Thus, $\nu\in \cM$, yielding closeness of $\cM$. 

\sp\fr So, since the set $\cM$ is closed and tight, its compactness follows from Prohorov's Compactness Theorem (Theorem 4.4 in \cite{Cra02}). The proof is complete. 
\epf

The last step in the proof of Theorem~\ref{thm 1} is the following.

\bprop\label{c.1}
Let $\cM$ be the invariant set of random measures coming from Lemma \ref{3.2}.
Then the map $\Phi : \cM \to \cM$ is continuous with respect to the narrow topology.
\eprop

\bpf
Suppose that $\Lambda$ is a directed set and $(\nu^\a)_{\a\in\La}$ is a net in $\cP_m (\Om)$ converging to some measure  $\nu \in \cM$ in the narrow topology. If 
$h_{\shift (x),\a} := 1 / \nu^\a _{\shift (x)} (\cL_x \1 )$, then
$$ 
\Phi(\nu^\a)_{x}
= \cL^*_x\left( \frac{1}{\nu^\a_{\shift (x)} (\cL_x\1 )} \nu^\a_{\shift (x)} \right)
=\cL^*_x\left( h_{\shift (x),\a}\;  \nu^\a_{\shift (x)} \right)\,.
$$ 
We have to estimate $h_{\shift (x),\a} $. Since $K_{0,\th (x)}\subset [F]_x$ and since $\nu^\al \in \cM$, we have
\begin{align*} 
\nu^\al_{\th (x)} (\pf_x\1 )
&\geq \int _{K_{0,\th (x)}} \pf_\om \1 d\nu^\al_{\th (x)} 
\geq \inf\lt(\(\pf_x \1\)\big|_{[F]_x}\rt) \; \nu^\al_{\th (x)}( K_{0,\th (x)})   
\geq \frac{1}{2c_F}
\end{align*} 
where $c_F^{-1}=\ess\!\inf\lt(\inf\lt(\(\pf_x\1\)\big|_{[F]_x}\rt):x\in X\rt)>0$ by (A) since the set $F$ is finite. Therefore,
\beq\label{3.17}
1/M \leq h_{\shift (x),\a}  \leq 2c_F\;.
\eeq
The measures $h_{\a}\nu^\a$, $\a\in\La$, need not be random probability measures but, thanks to \eqref{3.17},
their fiber total values $h_{\shift (x),\a}\nu^\a_{\shift (x)}(E_x^\infty)$ are uniformly bounded away from zero and $\infty$.
Since $\cM$ is compact, it follows that $\{h_{\a}\nu^\a:\a\in\La\}$ is a tight family.
Let $\mu$ be an arbitrary accumulation point of this net.
It has been shown (as a matter of fact for sequences but the same argument works for all nets) in Lemma~2.9 of \cite{RoyUrb2011} that then 
$$
\mu = h \nu
$$
for some measurable function $h:X\to (0, \infty )$. Since the dual operator $\pf^*:\cC_b^*(\Om)\to\cC_b^*(\Om)$ is continuous, the measure $\pf^*\mu$ is an accumulation point of the net $\(\pf^*\(h_{\a}\nu^\a\):\a\in\La\)$. But all elements of this net belong to $\cP_m(\Om)$, whence $\pf^*\mu\in \cP_m(\Om)$. This means that
the disintegrations of this measure 
$$
\cL^*_x\mu_{\shift (x)} = h_{\shift (x)} \cL^*_x ( \nu _{\shift (x)} ), \ x\in X,
$$
are all (Borel) probability measures in respective spaces $E_x^\infty$. Therefore, 
$$
1= \cL^*_x \mu_{\shift (x)} (\1)
= h_{\shift (x)} \cL^*_x( \nu _{\shift (x)} )(\1),
$$
which implies that 
$$
\mu _{\shift (x)}  
= \frac{1}{\cL^*_x (\nu _{\shift (x)} )(\1)} \nu_{\shift (x)}\ , \quad x\in X \,.
$$
This means that $\pf^*\mu=\Phi(\nu)$. Since also $\(\pf^*\(h_{\a}\nu^\a\):\a\in\La\)=\(\Phi(\nu^\a):\a\in\La\)$, we thus conclude that the net $\(\Phi(\nu^\a):\a\in\La\)$ converges to $\Phi(\nu)$. The proof of continuity of $\Phi$ is complete.
\epf

\bpf[Proof of Theorem \ref{thm 1}]
Consider the vector space $C_b^*(\Om)$ but now 
endowed with the narrow (weak convergence) topology analogously as for $\cP_m(X)$. It is standard to see that $C_b^*(\Om)$ becomes then a locally convex topological vector space. Since, by Lemma~\ref{c.2},
$\cM$ is a non-empty convex compact subset of $C_b^*(\Om)$,
since, by Lemma~\ref{3.2}, $\cM$ is $\Phi$--invariant and since, by Proposition~\ref{c.1}, $\Phi$ continuous, the Schauder-Tichonov Fixed Point Theorem
applies and yields a fixed point $\nu$ of $\Phi$ in $\cM$. But being such a fixed point $\nu\in \cM$ precisely means being a random conformal measure.
Finally, $\l_x =\pf_x^*\nu_{\th (x)} (E_x^\infty)$, and thus $\| \log \l\|_\infty <\infty$ results from the estimate \eqref{3.17}.
\epf

Here and in the sequel $\nu\in \cM$ is a conformal measure obtained in 
Lemma~\eqref{3.2}. We will frequently need the following useful estimate.

\blem\label{4.1} 
Suppose that $\hat\sg:\Om\to\Om$ and $\ph:\Om\to\R$ are fine. Then, for every finite set $D\subset E$ and for every integer $n\geq 0$ there exists a constant $\b =\beta_{D,n}>0$ such that
$$
\nu_x([\om]_x)\geq \b
$$
for every integer $0\leq k\leq n$, every tuple $\om\in D_x^k$ and $m$-a.e. $x\in X$.
\elem

\bpf
It is enough to show this statement for $k=n$. We keep the same notation as in Lemma \ref{3.2} and its proof. Since $\nu \in \cM$, by  Lemma~\eqref{3.2} (1) and (3.2), we have $\nu_x(K_{0,x})\geq 1- \ga_0(x ) = 1/2$ for $m$-a.e $x\in X$. Along with item (3.1) of Lemma~\eqref{3.2}, this yields $\nu_x( [F]_x)\geq 1/2$ for $m$-a.e. $x\in X$. In consequence, for each such $x$ there exists $b_x\in F$ such that 
\beq\label{4.2}
\nu_x( [b_x]_x)\geq \frac1{2q}, \quad \text{where} \quad q=\#F\,.
\eeq
Since both sets $D$ and $F$ are finite, by the topological mixing property of the map $\hat\sg:\Om\to\Om$, there exists an integer $N\ge 0$ such that, for every $x\in X$, every $a\in D$ and every $b\in F$ can be connected by a word $\tau=\tau(a,b)\in E_{\th(x)}^N$, i.e. $axb\in E_x^{N+2}$. We will use this fact in what follows with $b=b_x$, the point defined just above, see \eqref{4.2}.

Let $\om\in D_x^n$ be a word of length $n+1$ over the alphabet $D$ and let $b:=b_{\th ^{n+N}(x)}$. Let also $\tau:=\tau(\om_n,b)\in E_{\th^{n+1}(x)}^N$.
Then, by conformality of $\nu$,
\begin{align*} 
\nu_x([\om]_x)
&\geq \nu([\om\tau b]_x)
\geq \l_x^{-(n+N)}\exp\left(\inf\(S_{n+N}\ph_{\big | [\om\tau b]_x}\)\right)
\nu_{\th^{n+N}(x)}[b_{\th ^{n+N}(x) }]\,.
\end{align*} 
Notice that
$$
Q_n:=\inf \;\left\{ \exp \left(S_{n+N} \ph_{\big | [\om\tau b]_x}\right): \;
x\in X, \, \om\in D_x^n  \right\}
$$
is positive because of \eqref{ddd} since $D$ is finite and $b_{\th ^{n+N}(\om) }\in F$ varies in a finite set too.  
Therefore,
$$
\nu_x([\om]_x)\geq \frac1{2q}\exp\(-(n+N)\|\log \, \l\|_\infty\) Q_n>0\, \;.
$$
The proof is complete.
\epf

\section{Uniform bounds and two-norm inequality}
In this section we establish uniform  bounds for the iterated normalized operators 
$$
\npf^n_x= \l_x^{-n} \pf^n _x \quad \text{where} \quad \l_x^n= \l_x ...\l_{\th^{n-1} (x)} \,.
$$
We will then conclude from these estimates the existence of a random invariant measure equivalent to the random conformal measure produced in the previous section. 
From now on, the map $\hat\sg:\Om\to\Om$ and the potential $\ph:\Om\to\R$ are assumed again to be fine. We start with the following.

\bprop\label{4.3}
There exists a constant $\hat M<+\infty$ such that 
$$
\| \npf^n_x\|_\infty =\| \npf^n_x\1\|_\infty \leq\hat M \quad \text{for every $\; n\geq 0\; $ and every $\; x\in X$.}
$$
\eprop

\bpf
The first equality is obvious. Conformality of $\nu$ and the distortion result, Lemma~\ref{2.12}, yield, for every $\om\in E_{\th^n (x)}^\infty$,
$$
1
=\int_{E_{\th^n (x)}^\infty} \npf^n_x \1d\nu _{\th^n (x)}
\geq \frac{1}{K} \npf^n_x\1 (\om) 
\nu _{\th^n (x)} \big([\om_0]_{\th^n (x ) )}\big)
$$
Given $l\geq 1$, let $\b:= \b_{\{0,1,\ld,l\},0}$ come from Lemma~\ref{4.1}. Then 
$\nu _{\th^n (x)} \big([\om_0]_{\th^n (x) )}\big)\geq \b$ if $\om_0\in \{0,...,l\}$. For such points $\om$ we therefore get
\beq\label{4.4}
\npf^n_\om \1 (\om ) \leq K/\b \,.\eeq
By the property (B) and since $\log \l \in L^\infty$, we can adjust (increasing if necessary) the integer $l$ such that 
\beq\label{4.5}
\npf\1 _x (\om)\leq 1 
\eeq
for $m$-a.e. $x\in X$ if $\om\in E_{\th (x)}^\infty$ and $\om_0>l$.
We shall prove by induction that
\beq\label{4.5a}
\|\npf_x^n\1\|_\infty \leq K/\b
\eeq
for all integers $n\ge 1$ and $m$-a.e. $x\in X$.
For $n=1$ this directly results from \eqref{4.4} and \eqref{4.5}. For the inductive step, assume that \eqref{4.5a} holds for some $n\ge 1$. For all points $\om\in E_{\th^{n+1}(x)}^\infty$ with $\om_0\leq l$ everything is fine due to \eqref{4.4}.
So, suppose that $\om_0> l$. Then \eqref{4.5} along with the inductive hypothesis yield,
$$ 
\npf^{n+1}_x\1 (\om) 
= \npf _{\th^n (x)}\left( \npf^n_x \1\right) (\om)
\leq (K/\b) \npf_{\th^n (x)}\1 (\om)\leq K/\b\,.
$$
The proof is complete by taking $\hat M:=K/\b$.
\epf

\fr Having this uniform bound and the distortion Lemma \ref{2.9}, 
we can now obtain the following two-norm inequality. For every $x\in X$ and $n\ge 1$ let
$$
\l_x^n:=\l_x\l_{\th(x)}\cdots\l_{\th^{n-1}(x)}
$$
and
$$
\l_x^{-n}:=(\l_x^n)^{-1}.
$$
\blem\label{prop two norm} 
There exists a constant $S<\infty$ such that 
$$
v_\a (\npf_x^n g)\leq S \Big( \|g\|_\infty +  e^{-\al n} v_\a (g)\Big)\,.
$$
for every $x\in X$, every integer $n\geq 0$ and every $g\in \cH_\al\(E_{x}^\infty\)$.
\elem

\bpf
Let $x\in x$, $n\geq 0$, $g\in \cH_\al\(E_{x}^\infty\)$ and let $\om, \tau\in E_{\th^n(x)}^\infty$ with $\om_0=\tau_0\in E$. Then,
$$
|\npf_x^n g (\om)-\npf_x^n g (\tau)|\leq  \Sg_1+\Sg_2,
$$
where
$$
\Sg_2:=\l_x^{-n}\!\!\!\sum_{\substack{\g\in E_x^n \\ \g\om_0\in E_x^{n+1} }} 
\!\!\!\exp\(  S_n\ph(x,\g\tau)\) |g(\g\om)-g(\g\tau)|
\leq \npf_x^n \1 (\tau) v_\al (g) e^{-\al n }d^\a(\om,\tau)^\al
$$
and 
\begin{equation*}
    \begin{split} 
\Sg_1
:=& \l_x^{-n}\sum_{\substack{\g\in E_x^n \\ \g\om_0\in E_x^{n+1}}} 
\Big|  \exp\(S_n\ph(x,\g\om)\) - \exp\(S_n\ph(x,\g\tau)\) \Big||g(\g\om)|\\
\leq &\|g\|_\infty \l_x^{-n}\sum_{\substack{\g\in E_x^n \\ \g\om_0\in E_x^{n+1}}}        
 \exp\(S_n\ph(x,\g\om)\)
  \left|1- \frac{\exp\(S_n\ph(x,\g\tau)\)}{\exp\(S_n\ph(x,\g\om)\)} \right| \\
\leq &\|g\|_\infty  \| \npf_x^n\1\|_\infty C V_\a(\ph)d(\om,\tau)^\al,
    \end{split} 
\end{equation*}
where the last inequality results from Lemma~\ref{2.9} and the fact that $\ph:\Om\to\R$ is essentially $\a$-H\"older. It suffices now to combine the above estimates of both terms $\Sg_1$ and $\Sg_2$ along with the uniform bound from Proposition~\ref{4.3}. 
\epf

As an immediate consequence of the preceding lemma and of Proposition \ref{4.3}, we get the following.

\bprop \label{4.6}
For every integer $n\geq 0 $ and $m$-a.e. $x\in X$, we have that
$$
\| \npf^n_x\|_\a \leq \hat M,
$$
where here $\hat M$ is the maximum of $S$ and the former $\hat M$.
\eprop

\sp



\section{Invariant positive cones, Bowen's contraction and \\ the spectral gap}
From the uniform bounds obtained in the previous section one can construct
immediately invariant densities and measures by bare hands. However, in order to get further, more sophisticated, properties, especially a spectral gap, additional  investigations are needed.  Since
the phase spaces $E_x^\infty$ are not compact, we are here in a situation similar to that of random iterations of transcendental functions considered in \cite{MyUrb2014}. Therefore, we will now introduce, apply and
develop analogues of these tools worked out in the context of transcendental situation. Particularly, appropriated invariant positive cones
and employment of a Bowen's type contraction. Doing this, the invariant densities, hence the invariant measures equivalent to conformal measures, will emerge in a sense for free.

\subsection{Invariant cones}
Consider the following cones:
\beq\label{l21.5}
\cC_x:= \left\{g:E_x^\infty\to [0,+\infty) \; :\;  \|g\|_\infty \leq \cA \int g d\nu_x <\infty\; \; \text{and} \;\; v_\a (g)\leq H \int g d\nu_x\right\}\,.
\eeq
\beq\label{l21.6}
\cC_{x,0} := \left\{g\in \cC_x \; :\;\; g\leq 2 \hat M\cA \left(\int g d\nu_x\right) \; \npf _{\shift^{-1} (x)}\1 \right\}\,.
\eeq
Since we are primarily interested in the projective features of these cones, it is convenient to use  the following slices
\beq\label{l21.6slice}
\La _x = \{ g\in \cC_x \; , \;\; \nu_x (g) =1\} \quad \text{and} \quad \La _{x,0} =\La_x\cap \cC_{x,0} 
\; \; , \;\; x\in X\,.
\eeq
The constant $\hat M$, which is here and in the sequel, still forms a common upper bound of Proposition \ref{4.3}  and Proposition \ref{4.6}. Both type of cones do depend on the constants $\a \leq 1$, $\cA>0$, $H>0$. Those must be chosen carefully in the sequel in order to obtain invariant cones.

We continue to write $C$ for the distortion constant appearing in Lemma~\ref{2.9} and 
Lemma~\ref{2.12}.
First of all, fix an integer $k\geq 1$ so large that
\beq\label{l21.8}
\frac12 +\big(2\hat M+4 \big)e^{-\a k }\leq 1\,.
\eeq
Then, property (B) and the fact that $\|\log \l \|_\infty <\infty$ yield the existence of an integer $l_0\geq 0$ such that
\beq\label{l21.10}
2\hat M \npf_x \1 (x) \leq 1 \quad \text {for all $\om\in E_{\th (x)}^\infty$ with $\om_0\geq l_0$}
\eeq
for $m$-a.e. $x\in X$.
Let $\hat l_0\ge 0$ be the integer associated to $l_0$ by Lemma \ref{l1}. Iterating this procedure define further
$$
l_1:=\hat l_0 \, , \;\; l_2:=\hat l_1\; \; ...\;\; l_k:=\hat l_{k-1}.
$$
Set
\beq\label{5.0.2}
\cK_x:= \left\{\om\in E_x^\infty \; : \;\; \om_j \leq l_j \  \text{ for all } \ j=0,1,\ld, k\right\}\,, \quad x\in X\,.
\eeq
Then, for every $\om\in E_{x}^\infty \setminus \cK_x$, we have that $\om_0>l_0$.

Now, if $D=\left\{0,..., \max _{0\leq j\leq k} l_j\right\}$ and 
$\b=\b_{D,k}$ is from Lemma~\ref{4.1}, define
\beq \label{l21.13}
\cA:= 2\max\{ 1, \hat M, \b^{-1} \} \quad \text{and} \quad H:= 2\hat M\cA+4\,.
\eeq
Notice that $\cA \geq 1$. This ensures that the constant function $\1 \in \cC_x$, $x\in X$. 
Let finally $N_0\geq 1$ be such that 
\beq \label{m22.08}
\hat MHe^{-\a N_0} \leq 1\,.
\eeq

\bprop \label{l21.7}
With the above choice of constants and  for every $n\geq N_0$, we have
$$
\npf_x^n \left( \cC_x\right) \subset \cC_{\shift ^n (x),0} \subset \cC_{\shift ^n (x)} \; , \;\; x\in X\,.
$$
\eprop

\bpf
Let $g\in \cC_x$. We may assume without loss of generality that $\int g d\nu_x =1$.  
Fix $n\geq N_0$ arbitrarily.
We will show that $\npf_x^n g \in \cC_{\shift ^n (x) , 0}$. Clearly  $\npf_x^n g\geq 0$.
From the two-norm type inequality in Lemma~\ref{prop two norm} and from the definition of the cones, we get that
\beq \label{l21.21}
v_\a (\npf_x^n g)\leq \hat M\left( \cA + e^{-\a n}H\right)\leq \hat M\cA+1\leq H,
\eeq
where the last two inequalities result from the choice of $N_0$, i.e. \eqref{m22.08}, and from the definition of $H$.
Since $\nu_{\th^n(x)} \left( \npf_x^n g\right) =\nu_x (g) =1$, the second condition for belonging to $ \cC_{\shift ^n (x)}$ is thus satisfied by the function $\npf_x^n g$.
Also, by Proposition~\ref{4.3}, 
\beq \label{l21.22}
\npf_x^n g 
= \npf_{\shift^{n-1}(x)}\left( \npf_x^{n-1} g\right)\leq \hat M\|g\|_\infty  \npf_{\shift^{n-1}(x)}\1
\leq \hat M\cA \npf_{\shift^{n-1}(x)}\1\,.
\eeq
Hence, in order to conclude that $\pf_x^ng\in \cC_{\shift ^n (x),0} $ it remains to estimate $\| \npf_x^n g\|_\infty$ from above by $\cA$.
Since we already have \eqref{l21.21}, we obtain, for every $\om\in  \cK_{\th^{n}(x)}$, the following:
\begin{equation*}\begin{split}
1
&=\int \npf_x^ng \, d\nu_{\shift^n (x)}
\geq \int_{[\om|_k]_{\th^n(x)}}\npf_x^n g\, d\nu_{\shift^n (x)} \\
&\geq \left( \npf_x^n g(\om)-He^{-\a k} \right)\nu_{\shift^n (x)} ([\om |_k]_{\th^n(x)})\\
&\geq \left( \npf_x^n g(\om)-He^{-\a k} \right) \b \,,
\end{split}\end{equation*}
where the last inequality was written due to Lemma \ref{4.1}. Therefore, 
$$ 
\npf_x^n g(\om)
\leq \b^{-1}+He^{-\a k} \leq  \cA \left(\frac12 + (2\hat M+4)e^{-\a k} \right)
\leq \cA,
$$
by  \eqref{l21.13} and \eqref{l21.8}. 
If, on the other hand, $\om \not\in \cK_{\th^{n}(x)}$, then $\om_0\geq l_0$
and it follows from  \eqref{l21.22} and \eqref{l21.10} that
$$
\npf_x^ng(\om )\leq \hat M\cA \npf_{\shift^{n-1}(x)}\1 (\om)\leq \cA\,.
$$
The proof is complete.
\epf

\subsection{Cone Contraction via Bowen's Lemma.}
This part is based on a slightly stronger version of \eqref{l21.10}. Increasing $l_0\ge 1$ if necessary we may assume that 
\beq\label{m22.01}
2\cA\hat M \npf_x \1 (\om) < 1 \quad \text {for all $\om\in E_{\th (x)}^\infty$ with $\om_0\geq l_0$}
\eeq
for $m$-a.e. $x\in X$. We shall prove the following.

\blem\label{m22.02}
For every $l\geq l_0$ there are $N=N_l\geq N_0$ and $a=a_l>0$ such that 
$$
\npf_x^N g_{\big| [0,...,l]_{\th^n(x)}}\geq a \quad \text{for every} \ \, g\in \La_{x ,0}\ \, and \;\; x\in X\,.
$$
\elem

\bpf
Let $g\in \Lambda_{x,0}$. Since $\int gd\nu_x =1$, we have that $\|g\|_\infty \geq 1$. On the other hand, the choice of $l_0$ implies that
$$
g\leq 2\hat M\cA \, \npf_{\th^{-1} (x)}\1 \leq 1 \quad \text{in} \quad \left\{\om\in E_x^\infty:
\; \om_0\geq l_0 \right\}.
$$
Thus, there exists $\hat\om\in [0,...,l_0]_x$ with $g(\hat\om)\geq 1$.
Let $q\geq 1$ be so large that $He^{-\a q} \leq 1/2$. Since our random shift $\hat\sg:\Om\to\Om$ is fine, hence topologically mixing,
there exists $N=N_l\geq N_0$ such that $\hat\sg^N\big([\hat\om|_q]_x\big)\supset [0,...,l_0]_{\th^N (x)}$.
So, fixing an arbitrary $\tau\in [0,...,l_0]_{\th^N (x)}$, there exists $\hat\tau\in [\hat \om|_q]_x\cap \hat\sg^{-N} (\tau)$. We therefore get that
\begin{equation*}
\begin{split}
\npf_x^N g(\tau) 
&\geq \l_x^{-N}e^{S_N \ph_x(\hat\tau)}g(\hat\tau).
\end{split}
\end{equation*}
Now, $g(\hat\tau)\geq g(\hat\om) -H e^{-\a q}\geq 1-(1/2)=1/2$ and $\l_x^{-N}e^{S_N \ph_x(\hat\tau)}$ is bounded away from zero because of Lemma~\ref{l2} and inequality $||\log\la||_\infty<+\infty$ which is a part of Theorem \ref{thm 1}.
This shows the existence of a required constant $a>0$ and the proof is complete.
\epf

Notice that there is no way to get a global version of Lemma~\ref{m22.02}, valid on the whole shift space $\Om$. This is why we have to work with the following truncated functions:
$$
\chi_{_{l,x}}:= \1|_{ [0,...,l]_x} \npf_{\th^{-1} (x)}\1
\ , \quad l\geq 0 \;\;\text{and}\;\; x\in X\,.
$$ 
Then, 
\beq\label{1rs4.6}
0\leq \chi_{_{l,x}} \leq \npf_{\th^{-1} (x)}\1\quad \text{and}
\quad v_\a (\chi_{_{l,x}})\leq v_\a (\npf_{\th^{-1} (x)}\1).
\eeq
Moreover, we may suppose that $l_0\geq 1$ is sufficiently large so that
$$
\chi_{_{l,x}} \in \cC_{x, 0} \quad \text{for all} \quad l\geq l_0 \;\; \text{and $m$-a.e.}\;\; x\in X\,.
$$
We will now obtain a version of Bowen's contraction lemma~\cite[Lemma 1.9]{Bow75}. In order to do so, we have to define one more constant: let $\eta>0$ be such that
\beq\label{m22.03}
0<\eta\leq \min \left\{ \frac13 , \; \frac1H , \; \frac12\frac{a}{\hat M}\right\} \,.
\eeq


\blem\label{m22.04}
For every $l\geq l_0$ and with $N=N_l$ given by Lemma~\ref{m22.02},
$$ 
\npfx^N g - \eta \chi_{_{l,\shift^N(x)}} \in \cC_{{\shift^N(x)},0} \quad \text{for every}\quad g\in \La_{x,0}\,.
$$
\elem

\bpf
Let $x\in X$, let $g\in \La_{x,0}$, and let $l\geq l_0$.
Lemma~\ref{m22.02} and Proposition~\ref{4.3} (applied with $n=1$) show that for $0<\eta <2/(a\hat M)$,
$$
\npfx^N g -\eta \chi_{_{l,\shift^N(x)}}\le \frac{a}{2} >0 \quad
\text{on} \quad [0,\ld,l]_{\shift^N(x)}\,.
$$
Set
\beq\label{m22.05}
g':=\frac{\npf_x^N g-\eta\chi_{_{{l,\shift^N(x)}}}}{1-\eta_{_{l,\shift^N(x)}}}\,, \;\; \;\text{where}\;\; \;
\eta_{_{{l,\shift^N(x)}}}
:= \eta \int_{E_{\shift^N(x)}^\infty}\chi_{_{l,\shift^N(x)}}d\nu_{\shift^N(x)} \,.
\eeq
Then $\int g' d\nu_{\shift^N(x)}  =1$ and $g'>0$. Hence, by Proposition~\ref{4.3}, and since $g\in\La_x$, we have,
$$
\begin{aligned}
(1-\eta_{_{l,\shift^N(x)}})\,g'
&\le \npf_x^N(|g|)+\eta \chi_{_{l,\shift^N(x)}}
\leq  \hat M\|g\|_\infty \npf_{\shift^{(N-1)}(x)}\1 + \eta  \npf_{\shift^{(N-1)}(x)}\1 \\
&\leq (\hat M\cA +\eta )   \npf_{\shift^{(N-1)}(x)}\1 \,.
\end{aligned}
$$
But as $0<\eta_{_{l,\shift^N(x)}}\leq \eta\leq 1/3$, and $\eta\le 1/H$, the inequality above along with the definition of $H$ yield
$g'\leq 2\hat M\cA   \npf_{\shift^{(N-1)}(x)}\1 $. In conclusion, the function $g'\in \La_{{\shift^N(x)},0}$ provided that we can show that $g'\in \cC_{{\shift^N(x)},0} $.

In order to estimate the $\a$-variation of $g'$ we use the two-norm type inequality, i. e. Lemma~\ref{prop two norm}:
$$
v_\a(g') 
\leq \frac{1}{1-\eta_{_{l,\shift^N(x)}}} \left(\hat M\|g\|_\infty+\hat Mv_\a(g)e^{-\a N}+ \eta v_\a(\chi_{_{l,\shift^N(x)}})\right)\,.
$$
Remember that $g, \chi_{_{l,\shift^N(x)}}\in \cC_x$, that $\eta \leq \min\{1/3,1/H\}$, and that we have \eqref{m22.08}. Therefore,
$$
v_\a(g') \leq 2(\hat M\cA+1+1)=2M\cA+4=H\,.
$$
It remains to estimate $\|g'\|_\infty$. If $\om\in[0,\ld,l_0]_{\shift^N(x)}$, then
$$
1=\int_{E_{\shift^N(x)}^\infty} g' d\nu_{\shift^N(x)}
\geq \int_{[\om|_k]_{\shift^N(x)}}g' d\nu_{\shift^N(x)} 
\geq (g'(\om)- He^{-\a k} ) \nu _{\shift^N(x)}([\om|_k]_{\shift^N(x)})\,.
$$
Using Lemma~\ref{4.1} and the choice of $k$ in \eqref{l21.8}, we thus obtain
$$ 
g'(\om)
\leq \b{-1} + He^{-\a k}
\leq \frac{\cA}{2} + (2\hat M\cA+4)e^{-\a k} 
\leq \cA\left(\frac12 +(2\hat M+4)e^{-\a k}\right) 
\leq \cA\,.
$$
If $\om\in E_{\shift^N(x)}^\infty\sms [0,\ld,l_0]_{\shift^N(x)} $, then $g'(\om)\leq 2\hat M\cA\npf_{\shift^{(N-1)}(x)}\1 (\om)\leq \cA$
by the choice of $l_0$ (see \eqref{m22.01}). Thus $||g'||_\infty\le \cA$. Consequently $g'\in\La_{{\shift^N(x)},0}$, implying that $g\in \cC_{{\shift^N(x)},0}$. The proof is complete.
\epf

Applying repeatedly Lemma~\ref{m22.04} we now shall prove the desired contraction of Perron-Frobenius operators.

\bprop\label{m22.09}
For every $\ep >0$ there exists $n_\ep \geq 1$ such that for every $n\geq n_\ep$ and
$m$--a.e. $x\in X$,
\beq\label{28.21}
\left\|\npf _x^{n} g_x -\npf _x^{n} h_x\right\|_\a \leq\ep  \quad \text{for all} 
\quad g_x,h_x\in \La_{x,0} \,.\eeq
\eprop

\bpf
Let $l\geq l_0$ and $N=N_l\geq N_0$ be the same as in Lemma~\ref{m22.04}.
Let $g=g_x^{(0)}\in \La _{x,0}$. With the notation of the previous proof, and in particular with
the numbers $\eta_{_{l,\shift^N(x)}}$ defined in \eqref{m22.05}, we get from
Lemma~\ref{m22.04} that
$$
\npfx ^N g 
= \eta \chi_{_{l,\shift^N(x)}} + (1- \eta_{_{l,\shift^N(x)}})g_{\shift ^N (x)}^{(1)}
$$
for some $g_{\shift ^N (x)}^{(1)}\in \La _{\shift ^N (x),0}$. Applying $\npf_{\shift ^N(x)} ^N$ to this equation and using once more Lemma \ref{m22.04}, gives
$$
\npfx ^{2N} g
=  \eta \npf_{\shift ^N(x)} ^N\chi_{_{l,\shift^N(x)}}
+ (1- \eta_{_{l,\shift^N(x)}})\eta \chi_{_{l,\shift^{2N}(x)}}
+ (1- \eta_{_{l,\shift^N(x)}})(1- \eta_{_{l,\shift^{2N}(x)}}) g_{{\shift^{2N}(x)}}^{(2)}
$$
for some $g_{\shift ^{2N} (x)}^{(2)}\in \La _{\shift ^{2N} (x),0}$. Inductively, it follows that for every $k\geq 1$
there is a function $g_{\shift ^{kN} (x)}^{(k)}\in \La _{\shift ^{kN} (x),0}$ such that
$$
\npfx ^{kN} g
= \eta \sum_{j=1}^k \left( \prod_{i=1}^{j-1} (1-\eta_{_{l,\shift^{iN}(x)}})\right)
\npf_{\shift^{jN}(x)}^{(k-j)N}\chi_{_{l,\shift^{jN}(x)}} +  \prod_{i=1}^{k} (1-  \eta_{_{l,\shift^{iN}(x)}})\,g_{\shift ^{kN} (x)}^{(k)}\,.
$$
Observe that the first of these two terms does not depend on $g$. Therefore, for every two $g,h\in \La_{x,0}$
there are $g_{\shift ^{kN} (x)}^{(k)}, h_{\shift ^{kN} (x)}^{(k)}\in \La _{\shift ^{kN} (x),0}$ such that
\begin{equation}\label{120140501}
\npfx ^{kN} g -\npfx ^{kN} h 
=  \prod_{i=1}^{k} (1-\eta_{_{l,\shift^{iN}(x)}})\,
\left( g_{\shift ^{kN} (x)}^{(k)} - h_{\shift ^{kN} (x)}^{(k)}\right)\,.
\end{equation}
But, by property (A),  
$$
Q_{l_0}:=\ess\inf\(\npf_x\1|_{[0,\ld,l_0]_x}, x\in X\)>0.
$$
Therefore, we have for all $l\geq l_0$ and all $x\in X$ that
$$
\eta_{_{l,x}} 
= \eta \int_{E_x^\infty}\chi_{_{l,x}} d\nu_x
\geq \eta \int_{[0,\ld,l_0]_x} \npf_{\shift^{-1}(x)}\1 d\nu_x
\ge \eta Q_{l_0}\nu_x([0,\ld,l_0]_x)
\ge \eta Q_{l_0}\b_{l_0}>0,
$$ 
where $\b_{l_0}:=\b_{D,l_0}>0$, comes from Lemma~\ref{4.1}. Thus, writing, $\tilde\eta:=\eta Q_{l_0}\b_{l_0}>0$, we have that
$$
1-\eta_{_{l,x}}\le 1-\tilde \eta<1.
$$
Along with \eqref{120140501}, this allows us to deduce the uniform bound of Proposition~\ref{m22.09}, with some $n_{1, \ep}\ge 1$ sufficiently large, for the supremum norm rather than the H\"older one. In fact, in what follows we want $n_{1, \ep}$ to witness $\ep/\max\{2,4\hat M\}$ rather than merely $\ep$.
In order to get the appropriate estimate for the $\a$--variation, we need once more the two-norm type inequality, i. e. Lemma~\ref{prop two norm}. 
Write $n=m+n_{2,\ep}+n_{1, \ep}$ with some integer $n_{2, \ep}\ge 0$ to be determined in a moment and some integer $m\ge 0$. Then for all $g,h\in \La_{x,0}$, we have
$$
\begin{aligned}
v_\a\left(\npf_x^n g -\npf_x^nh \right)
&= v_\a\left( \npf_{\shift ^{n_{1, \ep}}(x)}^{m+n_{2, \ep}} \left( \npfx ^{n_{1, \ep}}( g -h)\right)\right)\\
& \leq \hat M\left( \left \| \npfx ^{n_{1, \ep}}( g -h)\right\| _\infty + e^{-\a(m+n_{2, \ep)}}v_\a\left( \npfx ^{n_{1,\ep}}( g -h)\right)\right) \\
&\leq  \hat M\frac{\ep}{4\hat M} + e^{-\a n_{2,\ep}}H \\
&=\frac{\ep}4+ 2e^{-\a n_{2,\ep}}H,
\end{aligned}
$$
where the second summand in the second last inequality was written due to the fact that $\npfx ^{n_{1, \ep}} g, \npfx ^{n_{1, \ep}} h \in \La_{\shift ^{n_{1, \ep}}(x) , 0}$. It suffices now to choose the integer $n_{2, \ep}\ge 0$ sufficiently large so that $2He^{-\a n_{2,\ep}}\le \ep/4$. Then $v_\a\left(\npf_x^n g -\npf_x^nh \right)\le \ep/2$, and in consequence
$$
||\npf_x^n g -\npf_x^nh||_\a\le \frac{\ep}2+\frac{\ep}2=\ep.
$$
The proof is complete.
\epf

\bpf[Proof of Theorem~\ref{thm 2} {\rm(1)}] 
For every integer $k\ge 0$ set $\rho^{(k)}:= \npf ^{k} \1$, i. e. $\rho_x^{(k)}=\npf_{\th^{-k}(x)}^{k}\1$. First of all, Proposition~\ref{l21.7} implies that $\rho^{(k)}_x\in\La _{x,0}$ for every $k\geq N_0$. Given $\ep>0$ arbitrary, put $q:=\max\{N_0,n_{\ep/2}\}$, the latter coming from Proposition~\ref{m22.09}. Proposition~\ref{l21.7} also implies that $\npf_{\th^{-n}(x)}^{n-q}\1\in\La_{\th^{-q}(x),0}$ for every integer $n\ge q+N_0$. Hence Proposition~\ref{m22.09} applies to give
$$
\big\|\rho ^q_x - \rho ^n_x\big\|_\a
=\big\|\npf_{\th^{-q}(x)}^{q}\1-\npf_{\th^{-n}(x)}^{n}\1\big\|
=\big\|\npf_{\th^{-k}(x)}^{q}\1-\npf_{\th^{-q}(x)}^{k}\(\npf_{\th^{-n}(x)}^{n-q}\1\)\big\|
\leq \ep/2
$$
for all $x\in X$. Therefore, if $k, l\ge q+N_0$, then
$$
\big\|\rho^l_x - \rho ^k_x\big\|_\a
\le \big\|\rho ^l_x - \rho ^q_x\big\|_\a+\big\|\rho ^q_x - \rho ^k_x\big\|_\a
\le \frac{\ep}2 + \frac{\ep}2
=\e
$$
for all $x\in X$.
This shows that $(\rho^{(n)}_x)_{n=0}^\infty$ is a Cauchy sequence in $(\cH_\a (E_x^\infty) , \| .\|_\b)$ uniformly with respect to $x\in X$. Hence, this sequence has a limit, call it $\rho_x$, in $\cH_\a (E_x^\infty)$, and $\rho_x\in \La _{x,0}$. Since also
$$
\npf_x\rho^{(n)}_x
=\npf_{\th^{-(n+1)}(\th(x))}^{n+1}\1
=\rho^{(n+1)}_{\th(x)},
$$ 
and since $\npf_x:\cH_\a (E_x^\infty)\to \cH_\a (E_{\th(x)}^\infty)$ is a continuous operator, we conclude that $\npf\rho_x = \rho_{\th(x)}$ for $m$-a.e. $x\in X$.
Uniqueness of this function follows immediately from the contraction formula \eqref{28.21} of Proposition~\ref{m22.09}.
\epf

\sp\bpf[Proof of Theorem~\ref{thm 2} {\rm (2)}]
Since $\cA,H \ge 2$, we have that
$$
\big\{\1+h_x: \;\|h_x\|_\a < 1/4 \big\} \subset \cC_x,
$$
for all $x\in X$.
Let $g\in \cH_\a (E_x^\infty)$, $g\not\equiv 0$, be arbitrary. Then
$$
h:= \frac{g}{8\|g\|_\a} =(h+\1)- \1
$$
is a difference of two functions from $\cC_x$. If $\ep >0$ and if $n=n_\ep $ is given by Proposition~\ref{m22.09}, then applying this proposition we get:
\begin{equation*}\begin{split}
\Bigg\| \npfx^n \Bigg(h -&\left(\int h d\nu_x\right) \den_x\Bigg) \Bigg\|_\a
=\left\| \npfx^n(\1+h)-\npfx^n\(\nu_x(\1+h)\den_x)-\npfx^n\1+\npfx^n\den_x\right\|_\a \le \\
&\leq \left\| \npfx^n(\1+h) - \npfx^n\(\nu_x(\1+h)\den_x)\right\|_\a
     +\left\|\npfx^n\den_x-\npfx^n\1\right\|_\a \\
&\leq \ep  \int (\1+h) \, d\nu_x + \ep \\
&\leq \frac{17}8\ep.
\end{split}\end{equation*}
This shows that for every $\ep >0$ there exists $N=N_\ep$ such that
$$
\left\| \npfx^N \left(g -\left(\int g d\nu_x\right) \den_x\right) \right\|_\a 
\leq \ep \| g\| _\a 
$$
for every $g\in \cH_\a (E_x^\infty)$. In particular, for the function $g-\nu_x(g)\rho_x$. So, we get
$$
\left\| \npfx^N \left(g -\left(\int g d\nu_x\right) \den_x\right) \right\|_\a 
\leq \ep \left\|g -\left(\int g d\nu_x\right)\rho_x\right\|_\a.
$$
Fix $\ep:=1/2$ and let $N=N_{1/2}$. Write any integer $n\ge 0$ in a unique form as $n=kN+m$, where $k\ge 0$ and $m\in \{0,...,N-1\}$.  Then, using also Proposition~\ref{4.6}, for every $g\in \cH_\a (E_x^\infty)$, we have
\begin{equation*}\begin{split}
\left\| \npfx^n g -  \int g d\nu_x \den_{\shift^n(x)} \right\|_\a 
&= 
\left\| \npf_{\shift^{kN}(x)}^m \left( \npfx^{kN}\big( g -\int g d\nu_x \den_x \big) \right)\right\|_\a\\
&\leq \hat M\left(\frac12\right)^k \left\| g -\int g d\nu_x \den_x\right\|_\a\\
&\leq 2\hat M\left(\frac1{2^{1/N}}\right)^n \left( 1+\|\den_x\|_\a\right) \| g\|_\a \,.
 \end{split}\end{equation*}
This completes the proof of Theorem~\ref{thm 2} {\rm (2)}.
\epf

\vfill
\pagebreak


\end{document}